\documentclass[12pt]{amsart}

\usepackage{amsmath}
\usepackage{amssymb}
\usepackage{amscd}
\usepackage{color}
\usepackage{hyperref}
\usepackage{mathrsfs}
\usepackage{enumitem}

\usepackage{mathtools}
\usepackage{eucal}
\usepackage{xy}
\xyoption{all}
\usepackage{bm}
\usepackage[normalem]{ulem}

\newcommand{\nc}{\newcommand}

\renewcommand{\AA}{{\mathbb{A}}}
\nc{\CC}{{\mathbb{C}}}
\renewcommand{\P}{{\mathbb{P}}}
\nc{\QQ}{{\mathbb{Q}}}
\nc{\ZZ}{{\mathbb{Z}}}

\nc{\cA}{{\mathcal{A}}}
\nc{\cB}{{\mathcal{B}}}
\nc{\cC}{{\mathcal{C}}}
\nc{\cD}{{\mathcal{D}}}
\nc{\cE}{{\mathcal{E}}}
\nc{\cF}{{\mathcal{F}}}
\nc{\cG}{{\mathcal{G}}}
\nc{\cI}{{\mathcal{I}}}
\nc{\cJ}{{\mathcal{J}}}
\nc{\cK}{{\mathcal{K}}}
\nc{\cL}{{\mathcal{L}}}
\nc{\cM}{{\mathcal{M}}}
\nc{\cN}{{\mathcal{N}}}
\nc{\cO}{{\mathcal{O}}}
\nc{\cP}{{\mathcal{P}}}
\nc{\cQ}{{\mathcal{Q}}}
\nc{\cR}{{\mathcal{R}}}
\nc{\cS}{{\mathcal{S}}}
\nc{\cU}{{\mathcal{U}}}
\nc{\cX}{{\mathcal{X}}}
\nc{\cY}{{\mathcal{Y}}}

\nc{\rc}{{\mathrm{c}}}

\nc{\rA}{{\mathrm{A}}}
\nc{\rB}{{\mathrm{B}}}
\nc{\rC}{{\mathrm{C}}}
\nc{\rE}{{\mathrm{E}}}
\nc{\rF}{{\mathrm{F}}}
\nc{\rG}{{\mathrm{G}}}
\nc{\rH}{{\mathrm{H}}}
\nc{\rK}{{\mathrm{K}}}
\nc{\rM}{{\mathrm{M}}}
\nc{\rP}{{\mathrm{P}}}

\nc{\bC}{{\mathbf{C}}}
\nc{\bH}{{\mathbf{H}}}
\nc{\bJ}{{\mathbf{J}}}
\nc{\bL}{{\mathbf{L}}}
\nc{\bO}{{\mathbf{O}}}
\nc{\bR}{{\mathbf{R}}}
\nc{\bS}{{\mathbf{S}}}
\nc{\bT}{{\mathbf{T}}}

\nc{\bfw}{{\mathbf{w}}}

\newcommand{\bal}{\bm{\alpha}}
\newcommand{\bmu}{\bm{\mu}}
\newcommand{\bka}{\bm{\kappa}}

\nc{\fM}{{\mathfrak{M}}}
\nc{\fS}{{\mathfrak{S}}}

\nc{\tX}{{\tilde{X}}}

\nc{\tcD}{{\tilde{\cD}}}
\nc{\tcP}{{\tilde{\cP}}}
\nc{\tcX}{{\tilde{\cX}}}
\nc{\bcX}{{\bar{\cX}}}

\nc{\tf}{{\tilde{f}}}

\nc{\eps}{\varepsilon}
\nc{\kk}{{\Bbbk}}

\nc{\Db}{{\mathbf{D}}^{\mathrm{b}}}
\nc{\Dp}{{\mathbf{D}}^{\mathrm{perf}}}
\nc{\FS}{\mathbf{FS}}
\nc{\Fuk}{\mathbf{Fuk}}
\nc{\MF}{\mathbf{MF}}

\nc{\QH}{\mathrm{QH}}
\nc{\QS}{\mathrm{QS}}
\nc{\BQH}{\mathrm{BQH}}

\newcommand{\dg}{{\mathrm{d}}}
\newcommand{\g}{{\mathrm{g}}}

\nc{\fMcub}{\fM_{\mathrm{cub}}}
\nc{\fMcubnod}{\fM_{\mathrm{cub}}^{\mathrm{nod}}}

\newcommand{\even}{{\mathrm{even}}}
\newcommand{\odd}{{\mathrm{odd}}}
\newcommand{\can}{{\mathrm{can}}}
\newcommand{\prim}{{\mathrm{prim}}}

\nc{\rDs}{\mathrm{D}_{\mathrm{short}}}

\newcommand{\CP}{\mathbb{CP}}

\DeclareMathOperator{\Aut}{\mathrm{Aut}}
\DeclareMathOperator{\Crit}{\mathrm{Crit}}

\DeclareMathOperator{\Hom}{\mathrm{Hom}}
\DeclareMathOperator{\Ext}{\mathrm{Ext}}

\DeclareMathOperator{\HOH}{\mathrm{HH}}

\DeclareMathOperator{\Spec}{\mathrm{Spec}}
\DeclareMathOperator{\Bl}{\mathrm{Bl}}
\DeclareMathOperator{\Sing}{\mathrm{Sing}}
\DeclareMathOperator{\Pic}{\mathrm{Pic}}
\DeclareMathOperator{\Cl}{\mathrm{Cl}}

\DeclareMathOperator{\Ker}{\mathrm{Ker}}

\DeclareMathOperator{\Cone}{\mathrm{Cone}}

\DeclareMathOperator{\Gr}{\mathrm{Gr}}
\DeclareMathOperator{\OGr}{\mathrm{OGr}}

\DeclareMathOperator{\Fl}{\mathrm{Fl}}

\DeclareMathOperator{\Spin}{\mathrm{Spin}}

\DeclareMathOperator{\id}{\mathrm{id}}

\theoremstyle{plain}

\newtheorem{theorem}{Theorem}[section]
\newtheorem{conjecture}[theorem]{Conjecture}
\newtheorem{lemma}[theorem]{Lemma}
\newtheorem{proposition}[theorem]{Proposition}
\newtheorem{corollary}[theorem]{Corollary}

\theoremstyle{definition}

\newtheorem{definition}[theorem]{Definition}
\newtheorem{example}[theorem]{Example}

\theoremstyle{remark}

\newtheorem{remark}[theorem]{Remark}

\nc{\hpd}{\natural}
\nc{\tbJ}{\widetilde{\bJ}}
\nc{\barV}{\bar{V}}
\nc{\barQ}{\bar{Q}}
\nc{\length}{\mathrm{length}}
\nc{\xrightiso}{ \xrightarrow{\ \raisebox{-0.5ex}[0ex][0ex]{$\sim$}\ }}

\DeclareMathOperator{\usdim}{{\overline{Sdim}}}
\DeclareMathOperator{\lsdim}{{\underline{Sdim}}}
\DeclareMathOperator{\ind}{ind}

\begin{document}

\title{Semiorthogonal decompositions in families}
\author{Alexander Kuznetsov}
\thanks{This work was performed at the Steklov International Mathematical Center 
and supported by the Ministry of Science and Higher Education of the Russian Federation (agreement no.~075-15-2019-1614).}
\address{Algebraic Geometry Section, 
Steklov Mathematical Institute of Russian Academy of Sciences, 8 Gubkin str., Moscow 119991 Russia}
\email{akuznet@mi-ras.ru}

\begin{abstract}
We discuss recent developments  in the study of se\-mi\-or\-tho\-go\-nal decompositions of algebraic varieties
with an emphasis on their behaviour in families.
First, we overview new results concerning homological projective duality.
Then we introduce residual categories, discuss their relation to small quantum cohomology, 
and compute Serre dimensions of residual categories of complete intersections.
After that we define simultaneous resolutions of singularities and describe a construction
that works in particular for nodal degenerations of even-dimensional varieties.
Finally, we introduce the concept of absorption of singularities which works under appropriate assumptions
for nodal degenerations of odd-dimensional varieties.
\end{abstract}

\maketitle

\section{Introduction}

Semiorthogonal decompositions of derived categories of algebraic varieties 
were introduced into the realm of algebraic geometry at the turn of the millennium by Bondal and Orlov~\cite{BO:preprint,BO:ICM};
since then the theory of semiorthogonal decompositions has become one of its central topics. 
Some advances in this theory have been surveyed in~\cite{K14};
in this sequel paper we discuss the progress obtained after~2014.

An important point of view in algebraic geometry, explained by Grothendieck, 
is that geometry should be studied in a relative situation.
Thus, the central object of algebraic geometry is a morphism~$\cX \to B$, 
i.e., a \emph{family of schemes}~$\{\cX_b\}_{b \in B}$ parameterized by the points~$b \in B$ of a base scheme.
Translating this point of view into the context of semiorthogonal decompositions, 
we understand that we should study semiorthogonal decompositions of schemes~$\cX/B$, 
especially \emph{$B$-linear semiorthogonal decompositions}.

The first step in this direction has been made in~\cite{K06}, 
where the notions of~$B$-linear triangulated categories and semiorthogonal decompositions have been introduced:
an enhanced triangulated category~$\cD$ is \textsf{$B$-linear} 
if it is endowed with a monoidal action of the monoidal category~$\Dp(B)$ of perfect complexes on~$B$.
For instance if~$f \colon \cX \to B$ is a scheme over~$B$, 
the bounded derived category of coherent sheaves~$\Db(\cX)$ is $B$-linear 
(where~$\cF \in \Dp(B)$ acts on~$\Db(\cX)$ by~$\cG \mapsto \cG \otimes f^*\cF$ for~$\cG \in \Db(\cX)$)
and a semiorthogonal decomposition
\begin{equation}
\label{eq:dbcx-sod}
\Db(\cX) = \langle \cD_1, \dots, \cD_n \rangle
\end{equation} 
is \textsf{$B$-linear} if~$\cD_i \otimes f^*(\Dp(B)) \subset \cD_i$ for each~$i$.

The next major step in this direction has been performed in~\cite{K11},
where the notion of base change for~$B$-linear semiorthogonal decompositions has been introduced:
given a~$B$-linear semiorthogonal decomposition~\eqref{eq:dbcx-sod} and a  morphism~$B' \to B$, 
a $B'$-linear semiorthogonal decomposition
\begin{equation*}
\Db(\cX \times_B B') = \langle \cD_{1, B'}, \dots, \cD_{n, B'} \rangle
\end{equation*}
has been constructed (under appropriate technical assumptions).
In particular, for each point~\mbox{$b \in B$} of the base scheme~$B$, the base change categories~$\cD_{i,b}$ are defined.
This allows one to consider a $B$-linear category~$\cD_i$ as a \emph{family of triangulated categories}~$\{\cD_{i,b}\}_{b \in B}$ 
parameterized by the points~$b \in B$ of the base scheme~$B$.

In this survey we discuss several independent topics, 
all of which, however, correspond to various semiorthogonal decompositions defined in families.

In sections~\S\S\ref{sec:hpd}--\ref{sec:residual} we discuss some standard material, 
developing and deepening results from~\cite{K14}.
First of all, we describe the main advances in the theory of homological projective duality (HPD) (see~\cite{K07}, \cite[\S3]{K14}),
namely \emph{categorical joins} and \emph{categorical cones}.
These categorical constructions provide appropriate homological counterparts of classical constructions of projective geometry
and are compatible with HPD.
This theory itself relies on a ``noncommutative'' (or rather categorical) version of HPD that was developed by Alex Perry in~\cite{Perry},
so we start~\S\ref{sec:hpd} with a short survey of noncommutative HPD;
\S\ref{ss:nhpd} can also serve as an introduction to HPD. 

After that in~\S\ref{ss:joins}, we introduce the construction of categorical joins 
and explain in what sense it is compatible with HPD,
then in~\S\ref{ss:nonlinear-hpd} we state the nonlinear HPD theorem, 
a generalization of the fundamental theorem of HPD 
in which linear sections of varieties are replaced by arbitrary intersections,
and then in~\S\ref{ss:cones} we introduce the construction of categorical cones (a particular case of categorical joins), 
and combining it with HPD for smooth quadrics and the nonlinear HPD theorem, we deduce the quadratic HPD theorem.
As an application of these results we deduce the duality conjecture for Gushel--Mukai varieties.
Finally, in~\S\ref{ss:other-hpd} we list a number of important developments in HPD not covered in this survey.

In~\S\ref{sec:residual} we introduce \emph{residual categories}: given a semiorthogonal decomposition
\begin{equation}
\label{eq:intro-residual}
\Db(X) = \langle \cR, \cB, \cB(1), \dots, \cB(m-1) \rangle,
\end{equation}
where the line bundle~$\cO_X(1)$ is an $m$-th root of the anticanonical line bundle of a Fano variety~$X$, 
i.e.,~$\omega_X^{-1} \cong \cO_X(m)$, and~$\cB$ is an admissible subcategory of~$\Db(X)$,
the leftmost component~$\cR$ of~\eqref{eq:intro-residual} is called the \textsf{residual category}.
In~\S\ref{ss:residual-ms} we explain a mirror symmetry interpretation of residual categories,
which justifies conjectures relating the structure of the semiorthogonal decomposition~\eqref{eq:intro-residual} and residual category~$\cR$ 
to the small quantum cohomology ring of~$X$, stated in~\S\ref{ss:conjectures}.
Further, in~\S\ref{ss:residual-homogeneous} we specify the predictions of the above conjectures for some homogeneous varieties,
namely Grassmannians and adjoint and coadjoint homogeneous varieties; 
remarkably, in all these cases the residual categories have a combinatorial nature:
when nonzero they are generated by completely orthogonal exceptional collections or equivalent to derived categories of Dynkin quivers.
Finally, in~\S\S\ref{ss:residual-hypersurfaces}--\ref{ss:residual-ci} 
we discuss the residual categories of hypersurfaces and complete intersections in projective spaces.
In these cases the structure of residual categories is much more complicated.
In the case of hypersurfaces the residual categories have a \emph{fractional Calabi--Yau property},
and in the case of complete intersections they are fractional Calabi--Yau up to an explicit spherical twist;
using this we show that they provide interesting examples of categories with distinct \emph{upper} and \emph{lower Serre dimensions}.

In the last part of this survey (\S\S\ref{sec:scr}--\ref{sec:absorption}) 
we discuss two ways to find a nice categorical replacement for a \emph{degeneration of schemes}, i.e., 
for a flat proper morphism~$f \colon \cX \to B$ to a smooth pointed curve~$(B,o)$ such that
\begin{itemize}
\item 
the morphism~$f^o \coloneqq f\vert_{\cX^o} \colon \cX^o \to B^o$ is smooth, and
\item 
the central fiber~$X \coloneqq \cX_o$ of~$f$ is singular.
\end{itemize}
Here and below we denote 
\begin{equation}
\label{eq:b-notation}
B^o \coloneqq B \setminus \{o\},
\qquad 
\cX^o \coloneqq \cX \times_B B^o,
\qquad 
\cX_o \coloneqq \cX \times_B \{o\},
\end{equation}
so that we have the following commutative diagram with Cartesian squares:
\begin{equation}
\label{eq:diagram-family}
\vcenter{\xymatrix@C=3em{
\hbox to 0em{\hss$X = {}$} \cX_o \ar@{^{(}->}[r]^i \ar[d] &
\cX \ar[d]^f &
\cX^o \ar@{_{(}->}[l]_j \ar[d]^{f^o}
\\
\{o\} \ar@{^{(}->}[r] &
B &
B^o \ar@{_{(}->}[l]
}}
\end{equation} 
We show that in the case where~$f$ is projective
and the central fiber~$X = \cX_o$ has only ordinary double points as singularities and under appropriate assumptions, 
one can find 
\begin{itemize}
\item 
a smooth and proper $B$-linear ``modification''~$\cD$ of the derived category~$\Db(\cX)$ of the total space,~and
\item 
a smooth and proper ``modification''~$\cD_o$ of the derived category~$\Db(\cX_o)$ of the central fiber,
\end{itemize}
such that~$\cD_o$ is the base change of~$\cD$ along the embedding~$\{o\} \hookrightarrow B$.
The precise meaning of the word ``modification'' depends on the parity of~$\dim(X)$
and is explained in the following two theorems.
For simplicity, we assume that the central fiber has a single ordinary double point.

In the case where~$\dim(X)$ is even we construct a simultaneous categorical resolution of singularities of~$\cX$,
which is a special case of a categorical resolution of singularities introduced in~\cite{K08} and~\cite[\S4]{K14}.

\begin{theorem}[{Theorem~\ref{theorem:nodal-general}}]
\label{thm:intro-even}
Assume given a commutative diagram~\eqref{eq:diagram-family} with Cartesian squares, 
where~$f$ is a flat projective morphism to a smooth pointed curve 
such that~$f^o$ is smooth and the central fiber~$X$ has a single ordinary double point~\mbox{$x_o \in X$}.

If~$\dim(X)$ is even and~$\cX$ has an ordinary double point at~$x_o$, 
there is a smooth and proper~$B$-linear triangulated category~$\cD$ and a commutative diagram
\begin{equation*}
\xymatrix@C=5em{
\Dp(X) \ar[d]_{\pi_o^*} \ar@<.5ex>[r]^{i_*} &
\Dp(\cX) \ar[d]_{\pi^*} \ar[r]^{j^*} \ar@<.5ex>[l]^{i^*} &
\Dp(\cX^o) \ar@{=}[d] 
\\
\cD_o \ar[d]_{\pi_{o*}} \ar@<.5ex>[r]^{i_*} &
\cD \ar[d]_{\pi_*} \ar[r]^{j^*} \ar@<.5ex>[l]^{i^*} &
\cD_{B^o} \ar@{=}[d] 
\\
\Db(X) \ar@<.5ex>[r]^{i_*} &
\Db(\cX) \ar[r]^{j^*} \ar@<.5ex>[l]^{i^*} &
\Db(\cX^o)
}
\end{equation*}
where~$\cD_o$ and~$\cD_{B^o}$ are the base change categories of~$\cD$, 
the functor~$\pi^*$ is left adjoint to~$\pi_*$, $\pi_o^*$ is left adjoint to~$\pi_{o*}$, and 
\begin{equation*}
\pi_* \circ \pi^* \cong \id,
\qquad 
\pi_{o*} \circ \pi_o^* \cong \id.
\end{equation*}
In particular, $\cD$ provides a categorical resolution of singularities for~$\cX$ 
and~$\cD_o$ provides a categorical resolution of singularities for~$X$.
\end{theorem}

In fact, even more is true --- the functors~$\pi^*$ and~$\pi_o^*$ are also right adjoint functors of~$\pi_*$ and~$\pi_{o*}$,
and the categorical resolutions~$\cD$ and~$\cD_o$ of~$\cX$ and~$X$ are \emph{weakly crepant} in the sense of~\cite[Definition~3.4]{K08};
this follows easily from the construction and~\cite[Proposition~4.5]{K08}.

The construction in the case where~$\dim(X)$ is odd is in some sense opposite to the above.
Note that in this case the exceptional divisor~$E_o \subset \Bl_{x_o}(X)$ of the blowup~~$\Bl_{x_o}(X)$
is a smooth quadric of even dimension, hence it comes with two \emph{spinor bundles}.

\begin{theorem}[{Corollary~\ref{cor:nmnf3} and Remark~\ref{rem:higher-and-even}}]
\label{thm:intro-odd}
Assume given a commutative diagram~\eqref{eq:diagram-family} with Cartesian squares, 
where~$f$ is a flat projective morphism to a smooth pointed curve 
such that~$f^o$ is smooth and the central fiber~$X$ has a single ordinary double point~\mbox{$x_o \in X$}.

If~$\dim(X)$ is odd, $\cX$ is smooth at~$x_o$, 
and there is an exceptional vector bundle~$\cE$ on~$\Bl_{x_o}(X)$ 
such that the restriction~$\cE\vert_{E_o}$ to the smooth quadric~$E_o$ is isomorphic to a spinor bundle
then there is a smooth and proper~$B$-linear triangulated category~$\cD$ and a commutative diagram
\begin{equation*}
\xymatrix@C=5em{
\langle \cP, \cD_o \rangle \ar@{=}[d] \ar@<.5ex>[r]^{i_*} &
\langle \langle i_*\cP \rangle, \cD \rangle \ar@{=}[d] \ar[r]^{j^*} \ar@<.5ex>[l]^{i^*} &
\cD_{B^o} \ar@{=}[d] 
\\
\Db(X) \ar@<.5ex>[r]^{i_*} &
\Db(\cX) \ar[r]^{j^*} \ar@<.5ex>[l]^{i^*} &
\Db(\cX^o)
}
\end{equation*}
where~$\cP \subset \Db(X)$ is an admissible subcategory, 
the triangulated subcategory~$\langle i_*\cP \rangle \subset \Db(\cX)$ generated by the image of~$\cP$ under~$i_*$ is admissible,
the categories~$\cD_o$ and~$\cD_{B^o}$ are base changes of~$\cD$,
the functors~$i_*$ and~$i^*$ in the top row are compatible with the semiorthogonal decompositions,
while~$j^*$ vanishes on~$\langle i_*\cP \rangle$. 
\end{theorem}

In both cases we obtain a smooth and proper $B$-linear category~$\cD$ 
such that for each point~$b \ne o$ in~$B$ the fiber~$\cD_b$ is equivalent to~$\Db(\cX_b)$, 
the derived category of the fiber of the original family of varieties.
Thus, the category~$\cD_o$ provides a \emph{smooth and proper extension}
of the family of categories~$\Db(\cX_b)$ across the point~$o \in B$.
Note however that this is achieved in two ``opposite'' ways --- 
in the situation described in Theorem~\ref{thm:intro-even} the category~$\cD_o$ is ``larger'' than~$\Db(X)$
(in particular, $\Db(X)$ is a Verdier quotient of~$\cD_o$),
while in the situation described in Theorem~\ref{thm:intro-odd} the category~$\cD_o$ is ``smaller'' than~$\Db(X)$
(in particular, $\cD_o$ is a semiorthogonal component of~$\Db(X)$).

This construction described in Theorem~\ref{thm:intro-odd}, which does not have a direct geometric analogue, 
is called \textsf{absorption of singularities}.
More precisely, we say that the subcategory~$\cP \subset \Db(X)$ \textsf{absorbs singularities} of~$X$ and, moreover,
it provides a \textsf{deformation absorption}.
We expect the notions of absorption and deformation absorption of singularities 
(defined in the general categorical context in~\S\ref{sec:absorption})
to become as important as the notion of resolution of singularities for the geometry of schemes.

The category~$\cP$ constructed in Theorem~\ref{thm:intro-odd} has a particularly interesting structure.
It is generated by a single object~$\rP \in \Db(X)$ such that 
\begin{equation*}
\Ext^\bullet(\rP,\rP) \cong \kk[t],
\qquad 
\deg(t) = 2.
\end{equation*}
We call such objects~\textsf{$\CP^\infty$-objects}.
They can be considered as limiting versions of~$\CP^n$-objects of Huybrechts--Thomas~\cite{HT06},
but while the latter give rise to autoequivalences of categories containing them,
the former provide semiorthogonal decompositions with interesting properties.
In particular, if~$\rP \in \Db(X)$ is a~$\CP^\infty$-object and~\mbox{$f \colon \cX \to B$} is a \emph{smoothing} of~$X$,
the object~$i_*\rP \in \Db(\cX)$ is exceptional.

We finish~\S\ref{sec:scr} and~\S\ref{sec:absorption} by sample applications of Theorem~\ref{thm:intro-even} and Theorem~\ref{thm:intro-odd}
to geometry of cubic fourfolds (see~\S\ref{ss:scr-cubic}) and Fano threefolds (see~\S\ref{ss:fano}), respectively.
In particular, we show that the nontrivial components of the derived categories of Fano threefolds of index~$2$ and degree~$1 \le d \le 5$
(with a minor modification in the case~$d = 1$)
can be represented as smooth and proper limits of the nontrivial components 
of the derived categories of Fano threefolds of index~$1$ and genus~$2d + 2$.
This gives a corrected version of a conjecture from~\cite{K09}.

\subsection*{Other important results}

Of course, this survey could not cover all interesting results related to semiorthogonal decompositions,
so we take this opportunity to list here some important achievements not mentioned in the body of the paper.

In a contrast to dimensions~$3$ and less, Fano varieties of higher dimensions are not yet classified.
However, there are several lists of interesting Fano varieties (e.g., see~\cite{Kuchle});
and of course it is interesting to describe their derived categories,
especially when they are expected to have interesting semiorthogonal components.
Some results in this direction can be found in~\cite{K15:Ku,K16:Ku,BFM,BS21}.
There is also some progress extending results about derived categories 
known over an algebraically closed field to more general fields, \cite{AB18,BDM19,K22}.

An interesting general question, directly related to the subject of this survey, 
is if it is possible to extend a semiorthogonal decomposition of a special fiber~$\cX_o$ of a family~$\cX/B$
to a $B$-linear semiorthogonal decomposition.
In~\cite{BOR} a positive answer is given under the assumption that~$\cX_o$ is smooth and proper,
and after an \'etale base change.

An intriguing connection between \emph{L-equivalence} of smooth projective varieties
(recall that~$X_1$ is \textsf{L-equivalent} to~$X_2$ if the difference of classes~$[X_1] - [X_2]$
is annihilated in the Grothendieck ring of varieties by the class~$[\AA^d]$ of an affine space)
and their derived equivalence was discovered, see~\cite{KS18,Efimov:L} and references therein.

Two important general results proved recently are Orlov's gluing theorem~\cite{Orlov16}
and Efimov's embeddability theorem~\cite{E20:phantoms}.
The first says that any \emph{gluing} of derived categories of smooth projective varieties
can be realized as an admissible subcategory of another smooth projective variety.
The second gives a criterion for realizability of an enhanced triangulated category
as an admissible subcategory in a category generated by exceptional collection.
It shows in particular that any phantom category admits such a realization,
thus providing a negative answer to~\cite[Conjecture~2.10]{K14}.

Finally, one of the mostly rapidly developing related areas is the study of Bridgeland stability conditions on semiorthogonal components.
We refer to~\cite{BLMS17,BM21} for surveys of this area.

\subsection*{Conventions}

In this paper all schemes are of finite type over a base field~$\kk$ and all categories are~$\kk$-linear.
We write~$\cA = \langle \cA_1, \dots, \cA_n \rangle$ for a semiorthogonal decomposition of a category~$\cA$ with components~$\cA_i$.
We denote by~$\Db(X)$ the bounded derived category of coherent sheaves and by~$\Dp(X)$ the category of perfect complexes on a scheme~$X$.
All pushforward, pullback, and tensor product functors are derived, although we use underived notation for them.
When we consider enhancements, we usually mean enhancements by differential graded categories as in~\cite{KL},
but one can also use infinity categories as in~\cite{Perry}.
We often use the notions of smoothness and properness for enhanced triangulated categories.
Recall that a dg-enhanced triangulated category is \textsf{smooth} if the diagonal bimodule 
over the underlying differential graded category is perfect,
and \textsf{proper} if for all objects~$\cF_1$, $\cF_2$ of the category 
the graded vector space~$\Ext^\bullet(\cF_1,\cF_2)$ has finite total dimension.

\subsection*{Acknowledgments}

Most of results described in this survey are obtained in collaboration.
I am very grateful to all my coauthors, especially to Alex Perry, Evgeny Shinder, and Maxim Smirnov for their invaluable inputs.
I would also like to thank all my colleagues for providing inspiration, sharing important ideas, and answering numerous questions.

Finally, I would like to thank Pieter Belmans, Alex Perry, Evgeny Shinder, and Maxim Smirnov
for their comments on the preliminary version of this paper.


\section{New results in homological projective duality}
\label{sec:hpd}

Homological projective duality studies the family of hyperplane sections of a given projective variety.
It was the main subject of the survey~\cite{K14} 
(see also~\cite{Thomas18} for an alternative perspective).
In this section we review the main advances in HPD obtained after~2014.

\subsection{Noncommutative HPD}
\label{ss:nhpd}

It has already become standard to consider nice triangulated categories as derived categories of ``noncommutative varieties''.
From this point of view it was clear from the very beginning
that the operation of homological projective duality is very noncommutative in nature ---
a manifestation of this is the fact that the result of HPD (even when applied to a commutative variety) is typically noncommutative.
However, it took some time for a firm foundation~\cite{Perry} for noncommutative HPD to be developed.

The setup of \emph{noncommutative HPD} is the following.
Instead of a smooth proper variety endowed with a morphism to a projective space~$\P(V)$ and a Lefschetz decomposition,
one considers a smooth and proper \textsf{Lefschetz category~$(\cA,\cA_0)$ over~$\P(V)$}.
By definition this consists of
\begin{itemize}
\item 
a \textsf{$\P(V)$-linear category}~$\cA$ (in the sense explained in the Introduction), and
\item 
an admissible subcategory~$\cA_0 \subset \cA$ (called \textsf{the Lefschetz center} of~$\cA$),
\end{itemize}
such that~$\cA_0$ extends to \textsf{right} and \textsf{left Lefschetz} (semiorthogonal) \textsf{decompositions}
\begin{equation*}
\cA = \langle \cA_0, \cA_1(1), \dots, \cA_{m-1}(m-1) \rangle 
\qquad\text{and}\qquad
\cA = \langle \cA_{1-m}(1-m), \dots, \cA_{-1}(-1), \cA_0 \rangle,
\end{equation*}
respectively,
where
\begin{equation*}
\cA_{m-1} \subset \dots \subset \cA_1 \subset \cA_0 
\qquad\text{and}\qquad
\cA_{1-m} \subset \dots \subset \cA_{-1} \subset \cA_0 
\end{equation*}
are two chains of admissible subcategories (called \textsf{the Lefschetz components} of~$\cA$).
The components~$\cA_i$ of both Lefschetz decompositions (if they exist) are determined by~$\cA_0$,
and if one of the Lefschetz decompositions exists then the other exists as well, \cite[Lemma~2.18, 2.19]{K08} or~\cite[Lemma~6.3]{Perry}.
Moreover, the maximal~$m$ such that~$\cA_{m-1} \ne 0$ equals the maximal~$m$ such that~$\cA_{1-m} \ne 0$; 
it is called the \textsf{length} of the Lefschetz category and is denoted~$\length(\cA)$.

The length of any Lefschetz category~$(\cA,\cA_0)$ over~$\P(V)$ 
satisfies the inequality
\begin{equation*}
\length(\cA) \le \dim(V),
\end{equation*}
and if the equality holds and if~$m = \length(\cA)$, the category~$\cA$ contains
\begin{equation}
\label{eq:rectangular-part}
\cA_{m-1} \otimes \Dp(\P(V)) = \langle \cA_{m-1}, \cA_{m-1}(1), \dots, \cA_{m-1}(m-1) \rangle
\end{equation}
as a rectangular Lefschetz subcategory (see~\cite[Corollary~6.19]{Perry})
and HPD for~$\cA$ reduces to HPD for the orthogonal of~\eqref{eq:rectangular-part} in~$\cA$ 
(this is the \emph{residual category} in the sense of~\S\ref{sec:residual}).
Thus, without losing generality, one can always reduce HPD to the case
where the Lefschetz category~$(\cA,\cA_0)$ is \textsf{moderate}, i.e., $\length(\cA) < \dim(V)$.

Given a Lefschetz category~$(\cA,\cA_0)$ over~$\P(V)$, one constructs the HPD Lefschetz category~$(\cA^\hpd,\cA^\hpd_0)$
by adapting the definition~\cite[Definition~6.1]{K07}.
Namely, consider the embedding of the universal hyperplane 
\begin{equation*}
\xymatrix@1{\bH(\P(V))\ \ \ar@{^{(}->}[r]^-\delta & \ \P(V) \times \P(V^\vee)}
\end{equation*}
and the base change~$\cA_{\bH(\P(V))}$ of the~$\P(V)$-linear category~$\cA$ 
along the natural projection~\mbox{$\bH(\P(V)) \to \P(V)$}.
Then~$\cA^\hpd$ is defined (see~\cite[Definition~7.1]{Perry}) as 
\begin{equation*}
\cA^\hpd \coloneqq \{ \cF \in \cA_{\bH(\P(V))} \mid \delta_*\cF \in \cA_0 \boxtimes \Db(\P(V^\vee)) \}
\subset \cA_{\bH(\P(V))},
\end{equation*}
which is $\P(V^\vee)$-linear category with respect to the $\P(V^\vee)$-linear structure of~$\cA_{\bH(\P(V))}$ 
induced by the morphism~$\bH(\P(V)) \to \P(V^\vee)$. 
Furthermore, the Lefschetz center
\begin{equation*}
\cA^\hpd_0 \subset \cA^\hpd 
\end{equation*}
is defined (see~\cite[Lemma~7.3]{Perry} and~\cite[(2.17)]{KP:joins}) 
as an explicit admissible subcategory in~$\cA^\hpd$.

For a linear subspace~$L \subset V$ we denote by
\begin{equation*}
L^\perp \coloneqq \Ker(V^\vee \to L^\vee) \subset V^\vee
\end{equation*}
its orthogonal subspace and 
we denote by~$\cA_{\P(L)}$ and~$\cA^\hpd_{\P(L^\perp)}$ the base change of~$\cA$ and~$\cA^\hpd$ 
along the embeddings~$\P(L) \hookrightarrow \P(V)$ and~$\P(L^\perp) \hookrightarrow \P(V^\vee)$, respectively.

The fundamental theorem of noncommutative HPD is stated as follows.

\begin{theorem}[{\cite[Theorem~8.7, 8.9]{Perry}}]
\label{thm:hpd-noncommutative}
Let~$(\cA,\cA_0)$ be a moderate Lefschetz category over~$\P(V)$.
Then the HPD Lefschetz category~$(\cA^\hpd,\cA^\hpd_0)$ over the dual projective space~$\P(V^\vee)$ is also moderate
and the HP double dual category is Lefschetz equivalent to the original 
\begin{equation*}
((\cA^\hpd)^\hpd, (\cA^\hpd)^\hpd_0) \simeq (\cA,\cA_0).
\end{equation*}
Moreover, if~$L \subset V$ is a linear subspace and~$L^\perp \subset V^\vee$ is its orthogonal with~$r = \dim(L)$ and~\mbox{$s = \dim(L^\perp)$},
and if~$m = \length(\cA)$, $n = \length(\cA^\hpd)$,
then there are semiorthogonal decompositions
\begin{align*}
\cA_{\P(L)} &= \langle \cK_L, \cA_s(1), \dots, \cA_{m-1}(m-s) \rangle,\\
\cA^\hpd_{\P(L^\perp)} &= \langle \cA^\hpd_{1-n}(r-n), \dots, \cA^\hpd_{-r}(-1), \cK'_{L^\perp} \rangle,
\end{align*}
and an equivalence of triangulated categories~$\cK_L \simeq \cK'_{L^\perp}$. 
\end{theorem}

The simplest example is \textsf{linear HPD} (see~\cite[\S8]{K07} for a relative version).

\begin{example}
\label{ex:hpd-linear}
Let~$0 \subsetneq W \subsetneq V$ be a linear subspace; thus~$\P(W)$ is a scheme over~$\P(V)$.
Then~$\Db(\P(W))$ endowed with the Lefschetz center~$\langle \cO_{\P(W)} \rangle \subset \Db(\P(W))$
is a moderate Lefschetz category over~$\P(V)$,
and the HPD of~$(\P(W),\langle \cO_{\P(W)} \rangle)$ is given by~$(\P(W^\perp),\langle \cO_{\P(W^\perp)} \rangle)$, 
where~$W^\perp \subset V^\vee$ is the orthogonal of~$W$.
\end{example}

See~\cite{K06,K06-gr,K07,K08-qu} for a number of other examples of HPD. 

Homological projective duality is related to classical projective duality via critical loci of morphisms~\cite[Theorem~7.9]{K07}
and this connection persists on the noncommutative level:
the classical projective dual of a Lefschetz category~$(\cA,\cA_0)$ 
(defined as the set of all hyperplanes in~$\P(V)$ such that the corresponding hyperplane section of~$\cA$ is singular)
coincides with the set of critical values of~$\cA^\hpd$ 
(defined as the set of points in~$\P(V^\vee)$ such that the corresponding fiber of~$\cA^\hpd$ is singular), \cite[Proposition~7.19]{Perry}.
When both~$\cA$ and~$\cA^\hpd$ are the derived categories of subvarieties~$X \subset \P(V)$ and~$Y \subset \P(V^\vee)$,
this reduces to classical projective duality~$X^\vee = Y$ (see Theorem~\ref{thm:hpd-quadrics} for an example).

Noncommutative HPD itself does not provide new examples of homologically projectively dual varieties (or categories)
but, as we already pointed above, it provides a firm background for developing the theory and for proving results 
like the one in the next subsection.

\subsection{Categorical joins}
\label{ss:joins}

The categorical join construction described below provides an appropriate homological extension 
of the classical join construction in projective geometry;
it is perfectly compatible with HPD; 
moreover, it provides new HPD examples and, as a consequence, 
new interesting results about derived categories of algebraic varieties.

Recall that the \textsf{join} of two projective varieties~$X_1 \subset \P(V_1)$ and~$X_2 \subset \P(V_2)$ 
is defined as the subvariety 
\begin{equation*}
\bJ(X_1,X_2) \subset \P(V_1 \oplus V_2)
\end{equation*}
swept out by all lines connecting points of~$X_1$ to points of~$X_2$, 
where we consider both~$X_1$ and~$X_2$ as subvarieties of~$\P(V_1 \oplus V_2)$ via the natural embeddings~$\P(V_i) \subset \P(V_1 \oplus V_2)$.
It is a well-known result in projective geometry that the join construction commutes with projective duality:
\begin{equation*}
\bJ(X_1,X_2)^\vee = \bJ(X_1^\vee,X_2^\vee) \subset \P(V_1^\vee \oplus V_2^\vee).
\end{equation*}
In~\cite{KP:joins} we define the categorical join of Lefschetz categories~$(\cA^1,\cA^1_0)$ over~$\P(V_1)$ 
and~$(\cA^2,\cA^2_0)$ over~$\P(V_2)$,
and prove a similar duality relation on the HPD level, see Theorem~\ref{thm:hpd-joins} below.
The definition is carried out in three steps.

In the \textbf{first step}, the join~$\bJ(X_1,X_2)$ is replaced by the \textsf{resolved join}
\begin{equation*}
\tbJ(X_1,X_2) \coloneqq \P_{X_1 \times X_2}(\cO(-1,0) \oplus \cO(0,-1)).
\end{equation*}
The resolved join is smooth (as soon as~$X_1$ and~$X_2$ are) 
and provides a natural resolution of singularities for the join~$\bJ(X_1,X_2)$, which is typically very singular.
In particular, we have the \textsf{universal resolved join}
\begin{equation*}
\tbJ(\P(V_1),\P(V_2)) = \P_{\P(V_1) \times \P(V_2)}(\cO(-1,0) \oplus \cO(0,-1)) \cong \Bl_{\P(V_1) \sqcup \P(V_2)}(\P(V_1 \oplus V_2)).
\end{equation*}
We denote by~$\eps_i \colon \P(V_1) \times \P(V_2) \hookrightarrow \tbJ(\P(V_1),\P(V_2))$ the exceptional divisor of the blowup lying over~$\P(V_i)$,
and by~$p \colon \tbJ(\P(V_1),\P(V_2)) \to \P(V_1) \times \P(V_2)$ the $\P^1$-bundle,
so that we have a commutative diagram
\begin{equation*}
\xymatrix@C=5em{
\P(V_1) \times \P(V_2) \ar[r]^{\eps_1} \ar[dr]_{\id} &
\tbJ(\P(V_1),\P(V_2)) \ar[d]^p &
\P(V_1) \times \P(V_2) \ar[l]_{\eps_2} \ar[dl]^{\id}
\\
&
\P(V_1) \times \P(V_2).
}
\end{equation*}
Note that the compositions~$p \circ \eps_i$ are isomorphisms.

In the \textbf{second step} we define the resolved join of~$\P(V_i)$-linear categories~$\cA^i$ as the base change
\begin{equation*}
\tbJ(\cA^1,\cA^2) \coloneqq (\cA^1 \boxtimes \cA^2)_{\tbJ(\P(V_1),\P(V_2))} 
\end{equation*}
of the~$\P(V_1) \times \P(V_2)$-linear category~$\cA^1 \boxtimes \cA^2$ along the~$\P^1$-bundle~$p$.
The blowup morphism
\begin{equation*}
\tbJ(\P(V_1),\P(V_2)) \to \P(V_1 \oplus V_2)
\end{equation*}
endows~$\tbJ(\cA^1,\cA^2)$ with a~$\P(V_1 \oplus V_2)$-linear structure.
The morphisms~$\eps_i$ and~$p$ defined above induce a commutative diagram of functors
\begin{equation*}
\xymatrix@C=5em{
\cA^1 \boxtimes \cA^2 &
\tbJ(\cA^1,\cA^2) \ar[l]_{\eps_1^*} \ar[r]^{\eps_2^*} &
\cA^1 \boxtimes \cA^2 
\\
&
\cA^1 \boxtimes \cA^2. \ar[ul]^{\id} \ar[ur]_{\id} \ar[u]_{p^*}
}
\end{equation*}
Note that the compositions~$\eps_i^* \circ p^*$ are equivalences.

So far, the construction uses the $\P(V_i)$-linear structure of the categories~$\cA^i$, 
but is independent of their Lefschetz centers~$\cA^1_0 \subset \cA^1$ and~$\cA^2_0 \subset \cA^2$;
they come into play in the \textbf{third step} of the construction.
We define the subcategories of~$\tbJ(\cA^1,\cA^2)$:
\begin{align*}
\cJ(\cA^1,\cA^2) &\coloneqq \{ \cF \in \tbJ(\cA^1,\cA^2) \mid 
\eps_1^*(\cF) \in \cA^1 \boxtimes \cA^2_0
\quad\text{and}\quad 
\eps_2^*(\cF) \in \cA^1_0 \boxtimes \cA^2
\},
\qquad\text{and}
\\
\cJ_0 &\coloneqq p^*(\cA^1_0 \boxtimes \cA^2_0).
\end{align*}
The isomorphisms~$\eps_i^* \circ p^* \cong \id$ imply the inclusion~$\cJ_0 \subset \cJ(\cA^1,\cA^2)$
and one can prove that the subcategory~$\cJ_0$ is a Lefschetz center in~$\cJ(\cA^1,\cA^2)$.

\begin{theorem}[{\cite[Theorem~3.21]{KP:joins}}]
\label{thm:joins-lefschetz}
If~$(\cA^1,\cA^1_0)$ and~$(\cA^2,\cA^2_0)$ are Lefschetz categories over projective spaces~$\P(V_1)$ and~$\P(V_2)$ 
then~$(\cJ(\cA^1,\cA^2), \cJ_0)$ is a Lefschetz category over~\mbox{$\P(V_1 \oplus V_2)$}
of length~$\length(\cA^1) + \length(\cA^2)$.
\end{theorem}

The categorical join can be thought of as a categorical resolution of the usual join, see~\cite[Proposition~3.17 and Remark~3.18]{KP:joins}.
The most important property of the categorical join operation is stated in the following

\begin{theorem}[{\cite[Theorem~4.1]{KP:joins}}]
\label{thm:hpd-joins}
If~$(\cA^1,\cA^1_0)$ and~$(\cA^2,\cA^2_0)$ are moderate Lefschetz categories over projective spaces~$\P(V_1)$ and~$\P(V_2)$ 
then there is a Lefschetz equivalence
\begin{equation*}
\cJ(\cA^1,\cA^2)^\hpd \simeq \cJ((\cA^1)^\hpd, (\cA^2)^\hpd),
\end{equation*}
where both sides are considered with their natural Lefschetz structures over~$\P(V_1^\vee \oplus V_2^\vee)$.
\end{theorem}

Many (but not all, see~\cite[\S6.2]{KP:joins}) geometric applications of these results 
rely on a categorical version of the following simple observation about the usual join operation.

\subsection{Nonlinear HPD theorem}
\label{ss:nonlinear-hpd}

As before assume given a pair of projective varieties~\mbox{$X_1 \subset \P(V_1)$} and~\mbox{$X_2 \subset \P(V_2)$}, 
but now assume~$\dim(V_1) = \dim(V_2)$.
Assume also given linear isomorphisms
\begin{equation*}
\xi_1 \colon V_1 \xrightiso V
\qquad\text{and}\qquad 
\xi_2 \colon V_2 \xrightiso V,
\end{equation*}
so that we can consider both~$X_1$ and~$X_2$ as subvarieties in~$\P(V)$.
Let 
\begin{equation*}
L(\xi_1,\xi_2) \coloneqq \Ker(V_1 \oplus V_2 \xrightarrow{\ (\xi_1, -\xi_2)\ } V)
\end{equation*}
be the equalizer of~$\xi_1$ and~$\xi_2$.
Then it is easy to check that
\begin{equation*}
\bJ(X_1,X_2) \cap \P(L(\xi_1,\xi_2)) = X_1 \cap X_2,
\end{equation*}
where we consider both sides as subvarieties in~$\P(V)$ using the identifications~$\xi_i$ 
and the induced identification~$L(\xi_1,\xi_2) \cong V$.
Furthermore, $X_1 \cap X_2$ can be thought of as~$X_1 \times_{\P(V)} X_2$,
and since the fiber product of varieties is categorified 
by the tensor product of linear categories (see~\cite[\S2.3]{Perry} for a definition), 
the above isomorphism has a categorical generalization:

\begin{lemma}[{\cite[Lemma~5.1]{KP:joins}}]
If~$(\cA^1,\cA^1_0)$ and~$(\cA^2,\cA^2_0)$ are Lefschetz categories over projective spaces~$\P(V_1)$ and~$\P(V_2)$ 
and~$\xi_i \colon V_i \xrightiso V$, $i = 1,2$, are isomorphisms, there is an equivalence of categories
\begin{equation*}
\cJ(\cA^1,\cA^2)_{\P(L(\xi_1,\xi_2))} \simeq
\cA^1 \otimes_{\Db(\P(V))} \cA^2,
\end{equation*}
between the base change of~$\cJ(\cA^1,\cA^2)$ along the inclusion~$\P(L(\xi_1,\xi_2)) \hookrightarrow \P(V_1 \oplus V_2)$
and the tensor product of the~$\P(V)$-linear categories~$\cA^1$ and~$\cA^2$ over~$\Db(\P(V))$.
\end{lemma}

Combining this observation with Theorem~\ref{thm:hpd-noncommutative} and Theorem~\ref{thm:hpd-joins} we obtain 

\begin{theorem}[{\cite[Theorem~5.5]{KP:joins}}]
\label{thm:hpd-nonlinear}
Let~$(\cA^1,\cA^1_0)$ and~$(\cA^2,\cA^2_0)$ be moderate Lefschetz categories 
over projective spaces~$\P(V_1)$ and~$\P(V_2)$ of equal dimensions.
Let
\begin{equation*}
N = \dim(V_1) = \dim(V_2),
\qquad
m = \length(\cJ(\cA^1,\cA^2)),
\qquad 
n = \length(\cJ((\cA^1)^\hpd,(\cA^2)^\hpd)), 
\end{equation*}
and let~$\cJ_i$ and~$\cJ_i^\hpd$ denote the Lefschetz components of~$\cJ(\cA^1,\cA^2)$ and~$\cJ((\cA^1)^\hpd,(\cA^2)^\hpd)$.
For any isomorphisms~$\xi_i \colon V_i \xrightiso V$, $i = 1,2$, there are semiorthogonal decompositions
\begin{align*}
\cA^1 \otimes_{\Db(\P(V))} \cA^2 &= 
\Big\langle \cK_{\xi_1,\xi_2}, \cJ_N(1), \dots, \cJ_{m-1}(m-N) \Big\rangle,
\\
(\cA^1)^\hpd \otimes_{\Db(\P(V^\vee))} (\cA^2)^\hpd &= 
\Big\langle \cJ_{1-n}^\hpd(N-n), \dots, \cJ_{-N}^\hpd(-1), \cK'_{(\xi_1^\vee)^{-1},(\xi_2^\vee)^{-1}} \Big\rangle,
\end{align*}
and an equivalence of categories~$\cK_{\xi_1,\xi_2} \simeq \cK'_{(\xi_1^\vee)^{-1},(\xi_2^\vee)^{-1}}$,
where the~$\P(V^\vee)$-linear structures of the categories~$(\cA^i)^\hpd$ 
are induced by the isomorphisms~$(\xi_i^\vee)^{-1} \colon V_i^\vee \xrightiso V^\vee$ for~$i = 1,2$.
\end{theorem}

In the special case where one of the Lefschetz categories, say~$\cA^2$, 
is the derived category of a linear subspace~$\P(L) \subset \P(V_2) \cong \P(V)$ with its natural Lefschetz structure,
then, as it was explained in Example~\ref{ex:hpd-linear}, 
$(\cA^2)^\hpd$ is Lefschetz equivalent to the derived category of the orthogonal subspace~$\P(L^\perp) \subset \P(V_2^\vee) \cong \P(V^\vee)$,
and there are equivalences
\begin{align*}
\cA^1 \otimes_{\Db(\P(V))} \Db(\P(L)) &\simeq \cA^1_{\P(L)},
\\
(\cA^1)^\hpd \otimes_{\Db(\P(V^\vee))} \Db(\P(L^\perp)) &\simeq (\cA^1)^\hpd_{\P(L^\perp)},
\end{align*}
where the right hand sides are the base change categories.
In this case the statement of Theorem~\ref{thm:hpd-nonlinear} is equivalent to the statement of Theorem~\ref{thm:hpd-noncommutative};
therefore Theorem~\ref{thm:hpd-nonlinear} can be considered as a \emph{nonlinear HPD theorem}.

Below we give two sample applications of Theorem~\ref{thm:hpd-nonlinear}.
The first is a consequence of HPD for the Grassmannian~$\Gr(2,5)$.

\begin{corollary}[{\cite[Proposition~1.1]{OR18}, \cite[Theorem~1.1]{BCP20}, \cite[Theorem~6.1]{KP:joins}}]
\label{cor:gr25-cy3}
Let~$V_1$ and~$V_2$ be vector spaces of dimension~$5$ and let~$\xi_i \colon \wedge^2V_i \xrightiso V$, $i = 1,2$, be linear isomorphisms.
Set 
\begin{equation*}
X \coloneqq \Gr(2,V_1) \times_{\P(V)} \Gr(2,V_2)
\qquad\text{and}\qquad 
Y \coloneqq \Gr(2,V_1^\vee) \times_{\P(V^\vee)} \Gr(2,V_2^\vee).
\end{equation*}
If the fiber products~$X$ and~$Y$ have the expected dimension~$3$, 
there is an equivalence of triangulated categories~\mbox{$\Db(X) \simeq \Db(Y)$}.
\end{corollary}

The varieties~$X$ and~$Y$ are deformation equivalent Calabi--Yau threefolds,
and as the above corollary states they are derived equivalent.
However, they are not birational in general~\cite[Theorem~1.2]{OR18}, \cite[Theorem~1.2]{BCP20},
and thus the pairs~$(X,Y)$ provide counterexamples to the so-called \emph{birational Torelli problem}. 
See~\cite{KaRa19} for another similar example.

The second application is a similar consequence of HPD for the connected components~$\OGr_\pm(5,10)$
of the orthogonal isotropic Grassmannian~$\OGr(5,10)$;
it provides examples of deformation and derived equivalent but not birational Calabi--Yau fivefolds.
Recall that~$\OGr_\pm(5,10)$ are homogeneous varieties of the simple algebraic group~$\Spin(10)$,
and the primitive generators of their Picard groups embed them into the projectivizations
of the two mutually dual $16$-dimensional half-spinor representations of~$\Spin(10)$.

\begin{corollary}[{\cite[Proposition~4.2]{Man19}, \cite[Theorem~6.3]{KP:joins}}]
\label{cor:ogr-cy}
Let~$V_1$ and~$V_2$ be vector spaces of dimension~$10$ endowed with non-degenerate quadratic forms,
let~$S_1$ and~$S_2$ be the~$16$-dimensional half-spinor representations of the corresponding groups~$\Spin(V_i)$,
and let~$\xi_i \colon S_i \xrightiso V$ be linear isomorphisms.
Set 
\begin{equation*}
X \coloneqq \OGr_+(5,V_1) \times_{\P(V)} \OGr_+(5,V_2)
\qquad\text{and}\qquad 
Y \coloneqq \OGr_-(5,V_1) \times_{\P(V^\vee)} \OGr_-(5,V_2).
\end{equation*}
If the fiber products~$X$ and~$Y$ have the expected dimension~$5$, 
there is an equivalence of triangulated categories~\mbox{$\Db(X) \simeq \Db(Y)$}.
\end{corollary}

Quadratic HPD discussed in~\S\ref{ss:cones} is also a special case of nonlinear HPD.

\subsection{Categorical cones and quadratic HPD}
\label{ss:cones}

As we have seen above, examples of geometrically meaningful Lefschetz categories 
for which the HPD categories are also geometrically meaningful and well understood
lead to applications of the nonlinear HPD theorem with interesting geometric consequences.
One nice example of such Lefschetz categories can be obtained from smooth quadrics.

Assume the base field~$\kk$ is algebraically closed of characteristic not equal to~$2$.
Let~$Q$ be a smooth quadric, i.e., a smooth proper variety isomorphic to a hypersurface of degree~$2$ in a projective space.
Let~$\cO_Q(1)$ be the ample line bundle that embeds~$Q$ as a quadric hypersurface.
A morphism~$f \colon Q \to \P(V)$ such that~$f^*\cO_{\P(V)}(1) \cong \cO_Q(1)$ is called \textsf{standard}:
thus either
\begin{itemize}
\item 
$f$ is a degree~$2$ embedding into a linear subspace of~$\P(V)$, or
\item 
$f$ is a degree~$2$ covering over a linear subspace of~$\P(V)$ ramified over a quadric.
\end{itemize}
We say that~$f$ is \textsf{non-degenerate} if the subspace above (i.e.\ the linear span of~$f(Q)$) is equal to~$\P(V)$.
In what follows we always consider~$Q$ as a~$\P(V)$-linear category by means of a standard morphism~$Q \to \P(V)$.

Let~$\cS$ be a spinor bundle on~$Q$ (the only one, if~$\dim(Q)$ is odd, or one of the two, if~$\dim(Q)$ is even).

\begin{lemma}[{\cite[Lemma~2.4]{KP21:quadrics}}]
\label{lemma:quadric-lefschetz}
The subcategory~$\cQ_0 \coloneqq \langle \cS, \cO \rangle \subset \Db(Q)$ is a Lefschetz center;
the length of the corresponding Lefschetz structure on~$\Db(Q)$ is equal to~$\dim(Q)$ and its Lefschetz components are given by
\begin{equation*}
\cQ_i = 
\begin{cases}
\langle \cS, \cO \rangle, & \text{if $|i| \le 1 - p$},\\
\langle \cO \rangle, & \text{if $1 - p < |i| \le \dim(Q) - 1$},
\end{cases}
\end{equation*}
where~$p \in \{0,1\}$ is the parity of~$\dim(Q)$.
\end{lemma}

The Lefschetz structure of~$\Db(Q)$ described above is called a \textsf{standard Lefschetz structure} of~$Q$;
note that it depends on the choice of the spinor bundle~$\cS$ 
(but the Lefschetz structures associated to different choices of~$\cS$ are noncanonically equivalent).

Recall that if~$Q \subset \P(V)$ is a smooth quadric hypersurface, the classical projective dual of~$Q$ 
is also a smooth quadric hypersurface~$Q^\vee \subset \P(V^\vee)$.
The HPD for a smooth quadric~$Q$ with a standard Lefschetz structure is described in similar terms.

\begin{theorem}[{\cite[Theorem~1.1]{KP21:quadrics}}]
\label{thm:hpd-quadrics}
Let~$f \colon Q \to \P(V)$ be a standard non-degenerate morphism of a smooth quadric~$Q$ endowed with a standard Lefschetz structure.
Then the HPD of~$Q$ is given by a standard non-degenerate morphism~\mbox{$f^\hpd \colon Q^\hpd \to \P(V^\vee)$}
of another smooth quadric~$Q^\hpd$, where:
\begin{enumerate}[label=\textup{(\arabic*)}]
\item 
\label{HPD-quadrics-1}
If~$f$ is a divisorial embedding and~$\dim(Q)$ is even, then~$Q^\hpd = Q^{\vee}$ is the classical projective dual of~$Q$
and~$f^\hpd \colon Q^{\hpd} \to \P(V^\vee)$ is its natural embedding. 
\item 
\label{HPD-quadrics-d-odd} 
If~$f$ is a divisorial embedding and~$\dim(Q)$ is odd, then~$f^\hpd \colon Q^{\hpd} \to \P(V^\vee)$ is the double covering
branched along the classical projective dual of~$Q$. 
\item 
\label{HPD-quadrics-dc-even}
If~$f$ is a double covering and~$\dim(Q)$ is even, then~$Q^\hpd$ is the classical projective dual of the branch locus of~$f$
and~$f^\hpd \colon Q^{\hpd} \to \P(V^\vee)$ is its natural embedding.
\item 
\label{HPD-quadrics-dc-odd} 
If~$f$ is a double covering and~$\dim(Q)$ is odd, then~$f^\hpd \colon Q^{\hpd} \to \P(V^\vee)$ is the double covering 
branched along the classical projective dual of the branch locus of~$f$. 
\end{enumerate} 
In all cases the HPD Lefschetz structure of~$Q^\hpd$ is a standard Lefschetz structure. 
\end{theorem}

This already allows one to apply the nonlinear HPD theorem, but the application becomes much more powerful
after an extension to singular quadrics (see~\cite[\S1.4]{KP:cones} for an explanation why this is useful).
Note that the derived category of a singular quadric is not smooth and proper, 
so it does not fit into the framework of HPD adopted in this paper.
On the other hand, every singular quadric~$Q$ can be written as a cone~$\bC_{V_0}(\barQ)$ over a smooth quadric~$\barQ$, 
and since a cone is a special case of a join,
one can use the formalism of categorical joins 
to find a suitable smooth and proper replacement for~$\Db(Q)$.
This is achieved by the \emph{categorical cone} construction.

Let~$(\cA,\cA_0)$ be a Lefschetz category over~$\P(V)$ and let~$V_0 \ne 0$ be a vector space.
We define the \textsf{categorical cone}~$\cC_{V_0}(\cA)$ as the categorical join
\begin{equation*}
\cC_{V_0}(\cA) \coloneqq \cJ(\P(V_0),\cA),
\end{equation*}
where~$\P(V_0)$ is endowed with the standard Lefschetz structure from Example~\ref{ex:hpd-linear}.
In fact, in~\cite{KP:cones} we use another definition, but it is equivalent to the above by~\cite[Proposition~3.15]{KP:cones}.

The categorical cone~$\cC_{V_0}(\cA)$, being the special case of a categorical join, 
is a Lefschetz category over~$\P(V_0 \oplus V)$, which is moderate if~$\cA$ is,
and can be thought of as a categorical resolution of the usual cone.
Combining Theorem~\ref{thm:hpd-joins} with Example~\ref{ex:hpd-linear} we deduce the following 

\begin{theorem}[{\cite[Theorem~1.1]{KP:cones}}]
\label{thm:hpd-cones}
Let~$V = V_0 \oplus \barV \oplus V_\infty$ and let~$(\cA,\cA_0)$ be a moderate Lefschetz category over~$\P(\barV)$.
Then there is a Lefschetz equivalence
\begin{equation*}
\cC_{V_0}(\cA)^\hpd \simeq \cC_{V_\infty^\vee}(\cA^\hpd),
\end{equation*}
where both sides are considered as Lefschetz categories over~$\P(V^\vee)$.
\end{theorem}

Here we add the summand~$V_\infty$ to~$V$ for higher flexibility of the construction.
For instance, in the next application to smooth quadrics it allows us to work  
with possibly degenerate morphisms of singular quadrics, like the morphism~$\bC_{V_0}(\barQ) \to \P(V)$ below.

\begin{corollary}[{\cite[Theorem~5.20]{KP:cones}}]
Let~$V = V_0 \oplus \barV \oplus V_\infty$, let~$\barQ \to \P(\barV)$ be a standard non-degenerate morphism of a smooth quadric~$\barQ$
and let~$\barQ^\hpd \to \P(\barV^\vee)$ be its HPD morphism.
There is a Lefschetz equivalence
\begin{equation*}
\cC_{V_0}(\barQ)^\hpd \simeq \cC_{V_\infty^\vee}(\barQ^\hpd),
\end{equation*}
where both sides are considered as Lefschetz categories over~$\P(V^\vee)$.

\end{corollary}

This leads to the following quadratic HPD theorem.
In the statement the assumption that the category~$\cA$ is supported away from~$\P(V_0)$ means that 
the~$\P(V)$-linear structure of~$\cA$ is induced by a~$(\P(V) \setminus \P(V_0))$-linear structure;
in this case the linear projection~\mbox{$\P(V) \setminus \P(V_0) \to \P(V/V_0)$} 
provides~$\cA$ with a~$\P(V/V_0)$-linear structure.
A similar convention is applied to~$\cA^\hpd$.

\begin{theorem}[{\cite[Theorem~5.21]{KP:cones}}]
\label{thm:hpd-quadratic}
Let~$(\cA,\cA_0)$ be a moderate Lefschetz category over~$\P(V)$ and let~$(\cA^\hpd,\cA^\hpd_0)$ be its HPD.
Assume given a direct sum decomposition~$V = V_0 \oplus \barV \oplus V_\infty$ 
such that the category~$\cA$ is supported away from~$\P(V_0)$ and~$\cA^\hpd$ is supported away from~$\P(V_\infty^\vee)$.
Let~$\barQ \to \P(\barV)$ be a standard non-degenerate morphism from a smooth quadric~$\barQ$ 
and let~$\barQ^\hpd \to \P(\barV^\vee)$ be its HPD morphism
and denote~$Q \coloneqq \bC_{V_0}(\barQ)$, $Q^\hpd \coloneqq \bC_{V_\infty^\vee}(\barQ^\hpd)$.
Let 
\begin{equation*}
N = \dim(V), \quad m = \length(\cA), \quad n = \length(\cA^\hpd), \quad 
d = \dim(Q), \quad e = \dim(Q^{\hpd}). 
\end{equation*} 
Then there are semiorthogonal decompositions
\begin{align*}
\begin{split}
\cA_{Q} = 
\Big \langle   \cK_Q(\cA) ,  
& \cA_{e}(1) \otimes \cS ,  \dots,  \cA_{m-1}(m-e) \otimes  \cS , \\
& \cA_{N-d}(1) \otimes  \cO ,  \dots,    \cA_{m-1}(m+d-N) \otimes  \cO   \Big\rangle.
\end{split}
\\
\begin{split}
(\cA^\hpd)_{Q^\hpd} = 
\Big\langle  
& \cA^\hpd_{1-n}(N-e-n) \otimes  \cO , \dots,   \cA^\hpd_{e-N}(-1) \otimes  \cO  ,\\   
& \cA^\hpd_{1-n}(d-n) \otimes  \cS^{\hpd} , \dots,  
\cA^\hpd_{-d}(-1) \otimes  \cS^{\hpd} ,  
\cK'_{Q^\hpd}(\cA^\hpd)   \Big\rangle , 
\end{split}
\end{align*}
and an equivalence of triangulated categories $\cK_Q(\cA) \simeq \cK'_{Q^\hpd}(\cA^\hpd)$,
where~$\cA_{Q}$ is defined as the base change of~$\cA$ 
along the morphism~$Q \to \P(V)$ and~$(\cA^\hpd)_{Q^\hpd}$ is defined analogously.
\end{theorem}

Again, here is a sample application of this result.
Recall that a \textsf{Gushel--Mukai variety}~\cite{DK1} is 
either a quadratic section of a linear section of~$\Gr(2,5)$, 
or a double covering of a linear section of~$\Gr(2,5)$ branched at a quadratic section.
In other words, a Gushel--Mukai variety can be described uniformly as a dimensionally transverse fiber product
\begin{equation*}
X = \Gr(2,V) \times_{\P(\wedge^2V)} Q,
\end{equation*}
where~$V$ is a $5$-dimensional vector space and~$Q \to \P(\wedge^2V)$ 
is a standard (possibly degenerate) morphism of a (possibly singular) quadric.
Note that for each~$Q$ as above there is a direct sum decomposition
\begin{equation}
\label{eq:wedge2v-decomposition}
\wedge^2V = V_0 \oplus \barV \oplus V_\infty
\end{equation}
and a standard non-degenerate morphism~$\barQ \to \P(\barV)$ from a smooth quadric~$\barQ$ 
such that one has~\mbox{$Q = \bC_{V_0}(\barQ)$}.

\begin{theorem}[{\cite[Theorem~6.4]{KP:cones}}]
Let~$V$ be a vector space of dimension~$5$, 
let~\eqref{eq:wedge2v-decomposition} be a direct sum decomposition of~$\wedge^2V$,
let~$\barQ \to \P(\barV)$ be a standard non-degenerate morphism of a smooth quadric, 
and let~$\barQ^\hpd \to \P(\barV^\vee)$ be its HPD morphism.
Assume the fiber products
\begin{equation}
\label{eq:gm-varieties}
X = \Gr(2, V) \times_{\P(\wedge^2V)} \bC_{V_0}(\barQ) \quad \text{and} \quad 
Y = \Gr(2, V^\vee) \times_{\P(\wedge^2 V^\vee)} \bC_{V_\infty^\vee}(\barQ^\hpd)
\end{equation}
are smooth GM varieties of dimensions~$d_X \geq 2$ and~$d_Y \geq 2$.
Let~$\cU_X$ and~$\cU'_Y$ denote the pullbacks to~$X$ and~$Y$ 
of the rank~$2$ tautological bundles of the corresponding Grassmannians.
Then there are semiorthogonal decompositions 
\begin{align}
\label{GMX}
\Db(X) & = \Big\langle \cK_X, \cO_X(1), \cU_X^{\vee}(1), \dots, \cO_X(d_X-2), \cU_X^{\vee}(d_X-2) \Big\rangle , \\ 
\label{GMY} 
\Db(Y) & = \Big\langle \cU'_Y(2-d_Y), \cO_Y(2-d_Y), \dots, \cU'_Y(-1), \cO_Y(-1), \cK'_Y \Big\rangle , 
\end{align} 
and an equivalence of triangulated categories $\cK_X \simeq \cK'_Y$. 
\end{theorem}

With a bit more work~\cite[Corollary~6.5]{KP:cones} 
this implies the duality conjecture~\cite[Conjecture~3.7]{KP:GM} for Gushel--Mukai varieties.

\subsection{Other results}
\label{ss:other-hpd}

To finish this section we list briefly other results developing HPD that appeared after~2014 and have not been mentioned in~\cite{K14}.
First, there are several works establishing HPD for new classes of varieties.
The most interesting among these are:
\begin{itemize}
\item 
The work of Rennemo~\cite{R20}, where the HPD for the symmetric square of a projective space~$\P^n$ (considered as a stack) is constructed,
see also~\cite{HT17} by Hosono--Takagi for a more geometric description of this HPD for small values of~$n$.
\item 
The work of Rennemo--Segal~\cite{RS19}, where a construction that allows to deduce some consequences of HPD for~$\Gr(2,2n+1)$
(without proving the HPD itself) is suggested.
\end{itemize}
Besides these major advances, the following papers should be mentioned:
\cite{BBF16}, where the linear HPD is applied to deduce HPD for determinantal varieties;
\cite{BDFIK}, where a differential graded algebra providing the HPD 
for a degree~$d$ (with~$d \ge 3$) Veronese embedding of a projective space is described;
and~\cite[\S D]{K19:dP6} where HPD for~$\P^1 \times \P^1 \times \P^1$ is established.

There are also some results contributing to general properties of HPD.
Among these one should mention~\cite{CT20}, where the HPD for a morphism~$f \colon X \to \P(V)$ 
is related to the HPD of the same variety~$X$ (blown up if necessary) with respect to a morphism~\mbox{$f' \colon X \to \P(V')$}
obtained from~$f$ by composing with a linear projection~$\P(V) \dashrightarrow \P(V')$
(see~\cite[\S B.1]{KP:joins} for a categorical version of this result).
Finally, one should mention the series of papers~\cite{Jiang21,Jiang:joins-intersections,Jiang:projectivizations,Jiang:blowups,Jiang:flops,Jiang:quot},
where an alternative approach to categorical joins is developed and many related results are obtained.

\section{Residual categories}
\label{sec:residual}

Let~$\cC$ be a triangulated category and let~$\bal_\cC \colon \cC \to \cC$ be an autoequivalence.
We say that an admissible subcategory~$\cB \subset \cC$ 
\textsf{generates a rectangular Lefschetz collection of length~$m$ with respect to~$\bal_\cC$}
if the collection of subcategories~$\cB, \bal_\cC(\cB), \dots, \bal_\cC^{m-1}(\cB)$ is semiorthogonal in~$\cC$. 
Admissibility of~$\cB$ implies that this collection extends to a semiorthogonal decomposition
\begin{equation}
\label{eq:residual}
\cC = \langle \cR, \cB, \bal_\cC(\cB), \dots, \bal_\cC^{m-1}(\cB) \rangle
\end{equation} 
and the component~$\cR$ of this decomposition is called the \textsf{residual category}, see~\cite{KS20,KP21}.

Residual categories often appear in HPD: if the Lefschetz decomposition of a Lefschetz category~$\cA$ 
is rectangular, i.e., $\cA_0 = \dots = \cA_{m-1}$,
the nontrivial components~$\cK_{\P(L)}$ appearing in Theorem~\ref{thm:hpd-noncommutative}, 
are the residual categories of~$\cA_{\P(L)}$.
In particular, residual categories often appear in families of semiorthogonal decompositions.

\subsection{Rotation functors}
\label{ss:rotation}

The residual category~$\cR$ defined by~\eqref{eq:residual} has especially nice properties 
if the subcategory~$\cB$ is \textsf{Serre compatible} in the sense that the condition
\begin{equation}
\label{eq:serre-compatibility}
\bS_\cC(\bal_\cC^m(\cB)) = \cB
\end{equation}
holds.
Note that the Serre functor commutes with any autoequivalence, hence~\eqref{eq:serre-compatibility} implies
that the autoequivalence~$\bS_\cC \circ \bal_\cC^m$ preserves all the components of~\eqref{eq:residual}.

One nice consequence of Serre compatibility is the following.
First of all, $\cR$ comes with a natural autoequivalence, induced by the so-called \textsf{rotation functor}~$\bO_\cB$.
Below we denote by~$\bL_\cB$ and~$\bR_\cB$ the left and right mutation functors of~$\cC$ with respect to~$\cB$.

\begin{proposition}[{\cite[Theorem~2.8]{KS20}, see also~\cite[Theorem~7.7]{KP17} and~\cite[Corollary~3.18]{K19}}]
If~$\cB \subset \cC$ is Serre compatible then the composition
\begin{equation*}
\bO_\cB \coloneqq \bL_{\cB} \circ \bal_\cC
\end{equation*} 
induces an autoequivalence of the residual category~$\cR$,
with the inverse autoequivalence induced by the composition~$\bal_\cC^{-1} \circ \bR_\cB$.
\end{proposition}

Second, the relation between the autoequivalence
\begin{equation*}
\bal_\cR \coloneqq \bO_\cB\vert_\cR
\end{equation*}
and the Serre functor~$\bS_\cR$ of~$\cR$ is analogous to that of~$\bal_\cC$ and~$\bS_\cC$.

\begin{theorem}[{\cite[Theorem~2.8, Remark~2.9, Proposition~2.10]{KS20}}]
\label{thm:residual-fine}
If~$\cB \subset \cC$ is Serre compatible and~$\cR$ is the residual category then
\begin{equation*}
\bS_\cR \circ \bal_\cR^m \cong (\bS_\cC \circ \bal_\cC^m)\vert_\cR.
\end{equation*}
Moreover, there is a bijection between 
\begin{itemize}
\item 
Lefschetz decompositions of~$\cR$ with respect to~$\bal_{\cR}$ and 
\item 
Lefschetz decompositions of~$\cC$ with respect to~$\bal_\cC$ containing~$\cB$ in every component.
\end{itemize}
\end{theorem}

\subsection{Mirror symmetry interpretation}
\label{ss:residual-ms}

Before discussing further properties of residual categories we sketch their interpretation from the point of view of mirror symmetry.
This section is mostly speculative, but it serves as a motivation for precise mathematical conjectures stated in~\S\ref{ss:conjectures}.

In this subsection we take~$\cC = \Db(X)$ to be the derived category of a smooth complex Fano variety~$X$ and let
\begin{equation}
\label{eq:balpha}
\bal_\cC(-) \coloneqq (-) \otimes \cL
\end{equation}
be the twist autoequivalence given by a line bundle~$\cL$.
Assume also that
\begin{equation}
\label{eq:omega-root}
\omega_X^{-1} \cong \cL^m,
\end{equation}
so that for any admissible subcategory~$\cB \subset \Db(X)$ generating a rectangular Lefschetz decomposition of length~$m$
the Serre compatibility condition~\eqref{eq:serre-compatibility} holds.

Homological mirror symmetry predicts the existence of a pair~$(Y,\bfw)$ 
(called a \textsf{Landau--Ginzburg model}) consisting of a proper morphism (called \textsf{the superpotential})
\begin{equation*}
\label{eq:lg-model}
\bfw \colon Y \to \AA^1
\end{equation*}
from a smooth scheme~$Y$ endowed with a relative symplectic form,
such that there are two equivalences of triangulated categories 
\begin{align}
\label{eq:db-fs}
\Db(X) &\simeq \FS(Y,\bfw),\\
\label{eq:fuk-mf}
\Fuk(X) &\simeq \MF(Y,\bfw),
\end{align}
where~$\Fuk(X)$ is the Fukaya category of~$X$, 
$\FS(Y,\bfw)$ is the Fukaya--Seidel category of~$(Y,\bfw)$, 
and~$\MF(Y,\bfw)$ is the category of matrix factorizations for~$(Y,\bfw)$.

A Landau--Ginzburg model for~$X$ is very far from being canonically defined.
On the other hand, the equivalence~\eqref{eq:db-fs} implies that the groups of autoequivalences of~$\Db(X)$ and~$\FS(Y,\bfw)$ coincide,
hence the symmetry of~$\Db(X)$ provided by the autoequivalence~\eqref{eq:balpha} 
should correspond to an autoequivalence~$\bal_\FS$ of~$\FS(Y,\bfw)$, i.e., to a symmetry of~$(Y,\bfw)$.
Since the (inverse) Serre functor of~$\FS(Y,\bfw)$ corresponds to the~$2\pi$-rotation around the origin 
of the target plane~$\AA^1 = \CC$ of the superpotential~$\bfw$,
the autoequivalence~$\bal_\FS$ should correspond to the~$2\pi/m$-rotation.
Therefore, we expect that there exists a \emph{$\bmu_m$-equivariant} Landau--Ginzburg model for~$X$,
i.e.\ a Landau--Ginzburg model~$(Y,\bfw)$, where~$Y$ is endowed with a~$\bmu_m$-action
and the morphism~$\bfw$ is~$\bmu_m$-equivariant for the standard~$\bmu_m$-action on~$\AA^1$.

So, from now on we assume that~$(Y,\bfw)$ is $\bmu_m$-equivariant.
The Fukaya--Seidel category~$\FS(Y,\bfw)$ is localized over the critical values of the superpotential~$\bfw$.
Let 
\begin{equation}
\label{eq:crit}
\Crit(\bfw) = \{z_0, z_1, \dots, z_N \} \subset \AA^1
\end{equation}
be the set of critical values of~$\bfw$.
For each~$0 \le i \le N$ choose a $\rC^\infty$-path~$\gamma_i$ in~$\AA^1$
connecting the point~$z_i$ to~$\infty = \P^1 \setminus \AA^1$
in such a way that the paths do not intersect 
and their natural cyclic order corresponding to the way they arrive at~$\infty$ 
is compatible with the linear ordering of the points~$z_i$ in~\eqref{eq:crit}.
By the definition of~$\FS(Y,\bfw)$ this gives a semiorthogonal decomposition 
(depending on the isotopy class of paths~$\gamma_i$)
\begin{equation*}
\FS(Y,\bfw) = \langle \cB_0, \cB_1, \dots, \cB_N \rangle
\end{equation*}
such that the component~$\cB_i$ is localized over~$z_i$.

Since~$\bfw$ is $\bmu_m$-equivariant, the set~\eqref{eq:crit} is~$\bmu_m$-invariant.
It may or may not contain the point~$0$;
in any case it will be convenient to set
\begin{equation*}
z_0 \coloneqq 0 \in \AA^1
\end{equation*}
and if this is not a critical value of~$\bfw$, set~$\cB_0 \coloneqq 0$.
So, let
\begin{equation*}
z_1,\dots,z_N \in \AA^1 \setminus \{0\} 
\end{equation*}
be the non-zero critical values of~$\bfw$.
Since the action of~$\bmu_m$ on~$\AA^1 \setminus \{0\}$ is free, $m$ divides~$N$ and reordering the points if necessary we can assume that
\begin{equation}
\label{eq:z-zeta}
z_{i + N/m} = \zeta \cdot z_i
\end{equation}
for~$1 \le i \le N - N/m$, where~$\zeta = \exp(2\pi\sqrt{-1}/m)$.
Furthermore, we can choose the paths~$\gamma_i$ in such a way that~$\gamma_{i + N/m} = \zeta \cdot \gamma_i$.
Then it follows that
\begin{equation*}
\cB_{i + N/m} = \bal_\FS(\cB_i),
\end{equation*}
so gathering the components~$\cB_1,\dots,\cB_{N/m}$ together and setting
\begin{equation*}
\cB \coloneqq \langle \cB_1, \dots, \cB_{N/m} \rangle,
\qquad 
\cR \coloneqq \cB_0
\end{equation*}
we see that~$\cB$ is admissible, Serre compatible of length~$m$,
and we have a semiorthogonal decomposition
\begin{equation*}
\FS(Y,\bfw) = \langle \cR, \cB, \bal_\FS(\cB), \dots, \bal_\FS^{m-1}(\cB) \rangle
\end{equation*}
with residual category~$\cR$.
In particular, if~$0$ is not a critical value for~$\bfw$, 
the residual category vanishes and~$\FS(Y,\bfw)$ acquires a rectangular Lefschetz decomposition.
In view of~\eqref{eq:db-fs}, the category~$\Db(X)$ should have a decomposition of the same type.

The above speculation shows the importance of understanding the critical values of the Landau--Ginzburg superpotential 
for the structure of~$\Db(X)$.
It is interesting that one can find them without describing the Landau--Ginzburg model, purely in terms of~$X$,
by using the second equivalence~\eqref{eq:fuk-mf}.
Indeed, the matrix factorization category by definition has a direct sum decomposition
\begin{equation*}
\MF(Y,\bfw) = \bigoplus_{i=0}^N \MF(Y,\bfw)_{z_i}
\end{equation*}
with components~$\MF(Y,\bfw)_{z_i}$ supported over the same points~$z_i \in \AA^1$,
therefore a similar direct sum decomposition holds for its Hochschild cohomology
\begin{equation*}
\HOH^\bullet(\MF(Y,\bfw)) = \bigoplus_{i=0}^N \HOH^\bullet(\MF(Y,\bfw))_{z_i}.
\end{equation*}
Thus, $\HOH^\bullet(\MF(Y,\bfw))$ can be thought of as a finite length coherent sheaf on~$\AA^1$,
i.e., a finite-dimensional module over the ring~$\CC[z]$ of functions on~$\AA^1$.
So, to control the points~$z_i$, it is enough to understand how the generator~$z$ of this ring acts on~$\HOH^\bullet(\MF(Y,\bfw))$.
For this we use the isomorphism
\begin{equation*}
\HOH^\bullet(\MF(Y,\bfw)) \cong \HOH^\bullet(\Fuk(X)) \cong \QH_\can(X)
\end{equation*}
(the first isomorphism follows from~\eqref{eq:fuk-mf} 
and the second has been conjectured in~\cite{Konts95} and is proved in~\cite[Corollary~9]{Ganatra}),
where~$\QH_\can(X)$ is the small quantum cohomology ring of~$X$ with the quantum parameters specialized to the anticanonical class
(see a more detailed discussion in the introduction to~\cite{KSm21}).
The right hand side~$\QH_\can(X)$ is isomorphic to~$\rH^\bullet(X,\CC)$ as a vector space
and is endowed with the supercommutative quantum multiplication;
in particular, one case consider the operator of quantum multiplication 
by the cohomology class~$\bka_X \in \QH_\can(X)$ of the anticanonical line bundle.
Note that cohomological degree (divided by~$2$) induces a~$\ZZ/m$-grading on the even part~$\QH_\can^\even(X)$ of~$\QH_\can(X)$
(which is a commutative ring), i.e., a~$\bmu_m$-action on its spectrum
\begin{equation}
\label{eq:qs}
\QS(X) \coloneqq \Spec(\QH_\can^\even(X)),
\end{equation}
and the corresponding morphism 
\begin{equation*}
\bka_X \colon \QS(X) \to \AA^1
\end{equation*}
is~$\bmu_m$-equivariant. 
It is expected that its image coincides with the~$\bmu_m$-invariant finite subset~\eqref{eq:crit} in~$\AA^1$,
see~\cite[Theorem~6.1]{Aur} for the toric case and a discussion preceeding it.

\subsection{The conjectures}
\label{ss:conjectures}

Summarizing the above discussion, we suggest the following precise conjectures.
We use the notation introduced in~\S\ref{ss:residual-ms}.
For a point~$z \in \AA^1$ we denote by~$\QH_\can^\even(X)_{\bka_X^{-1}(z)}$ the quotient ring of~$\QH_\can^\even(X)$
that corresponds to the union of connected components of~$\QS(X)$ supported over~$z$.

\begin{conjecture}
\label{conjecture:qh-db-general}
Let~$X$ be a complex Fano variety such that~\eqref{eq:omega-root} holds.
Assume the~$\bmu_m$-invariant subset
\begin{equation*}
\bka_X(\QS(X)) \cap (\AA^1 \setminus \{0\} ) = \{z_1,\dots,z_N\} \subset \AA^1 \setminus \{0\},
\end{equation*}
is ordered in such a way that~\eqref{eq:z-zeta} holds.
Then there is an $\Aut(X)$-invariant semiorthogonal decomposition 
\begin{equation}
\label{eq:residual-geometric}
\Db(X) = \langle \cR, \cB, \cB \otimes \cL, \dots, \cB \otimes \cL^{m-1} \rangle
\end{equation}
and the Hochschild homology of its components is given by
\begin{align*}
\HOH_\bullet(\cR) &= \hphantom{\ \bigoplus_{i=1}^{N/m}\ } \QH_\can(X) \otimes_{\QH_\can^\even(X)} \QH_\can^\even(X)_{\bka_X^{-1}(0)},\\
\HOH_\bullet(\cB) &= \ \bigoplus_{i=1}^{N/m}\ \QH_\can(X) \otimes_{\QH_\can^\even(X)} \QH_\can^\even(X)_{\bka_X^{-1}(z_i)},
\end{align*}
where in the right hand sides we identify~$\QH_\can(X)$ with the space~$\rH^\bullet(X,\CC) = \HOH_\bullet(\Db(X))$ 
with the Hochschild homology grading.
\end{conjecture}

The following example illustrates how this conjecture works.

\begin{example}
Let~$X \subset \P^n$ be a smooth Fano complete intersection of type~$(d_1,\dots,d_k)$ of dimension~\mbox{$\dim(X) \ge 3$}.
Then by~\cite[Theorem]{B95} the small quantum cohomology ring of~$X$ can be written as
\begin{equation*}
\QH_\can(X) = \CC\Big[h,\alpha\Big]_{\alpha \in \rH^{n-k}_\prim(X, \CC)} \Big/
\Big\langle h^{n-k+1} - h^\delta,\
h \cdot \alpha,\
\alpha_1 \cdot \alpha_2 - (\alpha_1,\alpha_2)(h^{n-k} - h^{\delta-1})
\Big\rangle,
\end{equation*}
where~$h$ is the hyperplane class, $\rH^{n-k}_\prim(X, \CC)$ is the primitive cohomology,
$(-,-)$ is the intersection pairing, 
and~$\delta = \sum (d_i - 1)$.
Since the anticanonical class of~$X$ is a multiple of~$h$, 
the localization of this ring away from~$\bka_X^{-1}(0)$ is obtained by inverting~$h$,
and so it is isomorphic to~$\CC[h]/(h^{n-k+1-\delta} - 1)$.
Note that~$m \coloneqq n-k+1-\delta = n + 1 - \sum d_i$ is the Fano index of~$X$, 
hence~$N = m$, $z_1,\dots,z_N$ are~$m$-th roots of unity, and~$\QH_\can(X)_{\bka_X^{-1}(z_i)} \cong \CC$ for~$1 \le i \le m$.
Thus, Conjecture~\ref{conjecture:qh-db-general} predicts the existence of semiorthogonal decomposition~\eqref{eq:residual-geometric}
with components~$\cB$ such that~$\HOH_\bullet(\cB) = \CC$.
Such a decomposition is indeed easy to construct, it is enough to take~$\cB = \langle \cO_X \rangle$, see~\S\ref{ss:residual-ci}.
\end{example}

It also makes sense to combine this conjecture with Dubrovin's conjecture~\cite{Dub}
that predicts that generic semisimplicity of the \emph{big quantum cohomology ring}~$\BQH(X)$
is equivalent to the existence of a full exceptional collection in~$\Db(X)$.
Note that generic semisimplicity of~$\BQH(X)$ implies that~$\rH^\odd(X,\CC) = 0$, 
hence~\mbox{$\QH_{\can}(X) = \QH^\even(X)$}, see~\cite[Theorem~1.3]{HMT}.
Recall the finite length scheme~$\QS(X)$ defined in~\eqref{eq:qs} and the~$\bmu_m$-equivariant morphism~$\bka_X$.
Furthermore, let
\begin{equation*}
\QS^\times(X) \coloneqq \bka_X^{-1}(\AA^1 \setminus \{0\}) \subset \QS(X),
\qquad
\QS^\circ(X) \coloneqq \QS(X) \setminus \QS^\times(X) \subset \QS(X).
\end{equation*}
These are finite~$\bmu_m$-invariant subschemes of~$\QS(X)$ and the action of~$\bmu_m$ on~$\QS^\times(X)$ is free.
Recall the autoequivalence~$\bal_\cR$ defined in~\S\ref{ss:rotation}.

\begin{conjecture}[{\cite[Conjecture~1.3]{KSm21}}]
\label{conj:qh-ec}
Let~$X$ be a complex Fano variety such that~\eqref{eq:omega-root} holds 
and the big quantum cohomology~$\BQH(X)$ is generically semisimple.
Let~$N$ be the length of the scheme~$\QS^\times(X)$.
\begin{enumerate}[label=\textup{(\roman*)}]
\item
There is a semiorthogonal decomposition~\eqref{eq:residual-geometric}, 
where the component~$\cB$ is generated by an~$\Aut(X)$-invariant exceptional collection of length~$N/m$.
\item
\label{conj:qh-ec:co}
The residual category~$\cR$ of~\eqref{eq:residual-geometric} has a completely orthogonal~$\Aut(X)$-invariant decomposition
\begin{equation*}
\cR = \bigoplus_{\xi \in \QS^\circ(X)} \cR_\xi
\end{equation*}
with components indexed by closed points $\xi \in \QS^\circ(X)$;
moreover, the component~$\cR_\xi$ of~$\cR$ is generated by an exceptional collection
of length equal to the length of the scheme~$\QS^\circ(X)$ at~$\xi$.
\item
The autoequivalence~$\bal_\cR$ permutes the components~$\cR_\xi$ of the residual category; 
more precisely, for each point~$\xi \in \QS^\circ(X)$ it induces an equivalence
\begin{equation*}
\bal_\cR \colon \cR_\xi \xrightiso \cR_{g(\xi)},
\end{equation*}
where~$g$ is a generator of~$\bmu_m$.
\end{enumerate}
\end{conjecture}

Most of the predictions in Conjecture~\ref{conj:qh-ec} are specializations of assertions of Conjecture~\ref{conjecture:qh-db-general}
to the case of a category with an exceptional collection.
The only exception is the complete orthogonality statement in part~\ref{conj:qh-ec:co}.
A justification for it, based on a comparison with the Fukaya--Seidel category~$\FS(Y,\bfw)$ 
can be found in~\cite{KS20} just before~\cite[Conjecture~1.12]{KS20}.

\begin{example}
Let~$X = \P^{m-1} \times \P^{m-1}$.
Then by the quantum K\"unneth formula~\cite{Kaufmann} one has 
\begin{align*}
\QH(X) &\cong \CC[h_1,h_2] / (h_1^m - q_1, h_2^m - q_2),
\\
\QH_\can(X) &\cong \CC[h_1,h_2] / (h_1^m - 1, h_2^m - 1),
\end{align*}
and the $\ZZ/m$-grading is defined by~$\deg(h_1) = \deg(h_2) = 1$.
Therefore,
\begin{equation*}
\QS(X) = \bmu_m \times \bmu_m \subset \AA^2,
\end{equation*}
where the embedding is induced by the natural embeddings~$\bmu_m \subset \AA^1 \setminus \{0\} \subset \AA^1$,
and up to rescaling of~$\AA^1$ the map~$\bka_X$ is induced by the summation map~$\AA^2 \to \AA^1$.

If~$m$ is odd, $\QS^\circ(X)$ is empty, 
hence Conjecture~\ref{conj:qh-ec} predicts the existence of an~$\Aut(X)$-invariant rectangular Lefschetz collection with zero residual category.
Several such collections, indeed, have been constructed in~\cite{R20}.

On the other hand, if~$m$ is even, $\QS^\circ(X)$ has length~$m$ (it consists of all pairs~$(\xi,-\xi)$ for~\mbox{$\xi \in \bmu_m$});
consequently, Conjecture~\ref{conj:qh-ec} predicts the existence of an~$\Aut(X)$-invariant rectangular Lefschetz collection 
with residual category generated by~$m$ completely orthogonal objects.
Such collection has been found in~\cite[Example~1.4]{KSm21}.
\end{example}

\subsection{Residual categories of homogeneous varieties}
\label{ss:residual-homogeneous}

In this subsection we make the predictions of Conjecture~\ref{conj:qh-ec} more explicit 
for homogeneous varieties of reductive algebraic groups
and compare them with known results about their derived categories.

According to an old folklore conjecture homogeneous varieties are expected to have full exceptional collections 
(the conjecture is still not proved, see a discussion of known cases in~\cite[\S1.2]{KP16} 
and more recent developments in~\cite{G20,F19,KSm21,BKS21,S21}).
Therefore, Conjecture~\ref{conj:qh-ec} should be applicable and we only need to compute the small quantum cohomology ring.
This is pretty easy for Grassmannians.

\begin{example}
\label{ex:qh-gr}
The small quantum cohomology ring of~$X = \Gr(k,n)$ can be presented as
\begin{equation*}
\QH(X) = \CC[c_1,\dots,c_k,s_1,\dots,s_{n-k}] / \langle (1 + c_1 + \dots + c_k)(1 + s_1 + \dots + s_{n-k}) = 1 + (-1)^k q \rangle,
\end{equation*}
where~$c_i$ and~$s_j$ should be thought of as Chern classes of the tautological subbundle and quotient bundle respectively,
and~$q$ is the quantum parameter (see~\cite[Theorem~0.1]{ST97}).
Decomposing formally
\begin{equation*}
1 + c_1 + \dots + c_k = \prod_{i=1}^k(1 - x_i),
\qquad 
1 + s_1 + \dots + s_{n-k} = \prod_{i=k + 1}^n(1 - x_i),
\end{equation*}
where~$x_1,\dots,x_n$ are the corresponding Chern roots, and specializing the quantum parameter~$q$ to~$(-1)^{k+1}$, 
we conclude that~$\{x_1,\dots,x_n\} = \bmu_n$.
Since the canonical class of~$X$ is proportional to the first Chern class of the tautological bundle,
we have up to rescaling 
\begin{equation*}
\QS(X) = \Big((\bmu_n)^k \cap ( \AA^k \setminus \Delta ) \Big) \Big/ \fS_k
\end{equation*}
where~$\Delta \subset \AA^k$ is the big diagonal, $\fS_k$ is the permutation group, 
and the map~$\bka_X$ is induced by summation of coordinates in~$\AA^k$.
Furthermore, the natural action of~$\bmu_n$ on~$\QS(X)$ is given 
by simultaneous multiplication of all coordinates by a root of unity.
Thus, the $\bmu_n$-action is free if and only if~$\gcd(k,n) = 1$, 
and otherwise orbits of length~$d$ correspond to subsets of cardinality~$k/d$ in~$\bmu_{n/d}$.

Note that some free orbits may be contained in~$\QS^\circ(X)$ 
(by~\cite{Gary} this happens if and only if both~$k$ and~$n - k$ are sums of non-trivial divisors of~$n$,
the simplest example with~$\gcd(k,n) = 1$ being~$n = 12$ and~$k = 5$),
so the following conjecture is stronger than the prediction of Conjecture~\ref{conj:qh-ec}.
\end{example}

\begin{conjecture}[{\cite[Conjecture~3.10 and Lemma~3.9]{KS20}}]
If~$X = \Gr(k,n)$ there is a rectangular Lefschetz collection with residual category generated by 
\begin{equation*}
R_{k,n} = - \sum_{d \,\vert \gcd(k,n),\ d > 1} \mu(d) \binom{n/d}{k/d}
\end{equation*}
completely orthogonal objects, where
\begin{equation*}
\mu(d) \coloneqq
\begin{cases}
\hphantom{-}1, & \text{if~$d$ is square-free with an even number of prime factors},\\
-1, & \text{if~$d$ is square-free with an odd number of prime factors},\\
\hphantom{-}0, & \text{if $d$ has a squared prime factor}
\end{cases}
\end{equation*}
is the M\"obius function.
\end{conjecture}

By now this conjecture is known for the case~$\gcd(k,n) = 1$, where the residual category vanishes~\cite[Theorem~4.3 and Proposition~4.8]{F13}, 
as well as for the case where~$k$ is a prime number~\cite[Theorem~3.13]{KS20}.

For more complicated homogeneous varieties we can also use the available quantum cohomology computations 
to formulate a number of precise conjectures.
We do this below for some interesting classes of homogeneous varieties.

Recall that the \textsf{adjoint} (resp.\ \textsf{coadjoint}) homogeneous variety of a simple algebraic group $\rG$
is the orbit of the highest weight vector in the projectivization of an irreducible representation of~$\rG$,
whose highest weight is the highest \emph{long} (resp.\ \emph{short}) root of~$\rG$;
in particular, if the group~$\rG$ is simply laced, the adjoint and coadjoint varieties coincide.

The following conjecture is motivated by the results of~\cite{PeSm}.
For a group~$\rG$ we denote by~$\rDs(\rG)$ the short roots subdiagram of the Dynkin diagram of~$\rG$
(if the group~$\rG$ is simply laced, it is the entire Dynkin diagram).

\begin{conjecture}[{\cite[Conjecture~1.8]{KSm21}}]
Let~$X$ be the coadjoint variety of a simple complex algebraic group~$\rG$.
Then $\Db(X)$ has an~$\Aut(X)$-invariant rectangular Lefschetz exceptional collection with residual category~$\cR$, where
\begin{enumerate}[label=\textup{(\arabic*\textup)}]
\item
if the Dynkin type of~$\rG$ is~$\rA_n$ and $n$ is even, then $\cR = 0$;
\item
otherwise, $\cR$ is equivalent to the derived category of representations of a quiver of type~$\rDs(\rG)$.
\end{enumerate}
\end{conjecture}

By now this conjecture is known for all Dynkin types except for the exceptional types~$\rE_6$, $\rE_7$, $\rE_8$, 
see a discussion and references in~\cite{KSm21}.

For adjoint varieties (of non simply laced groups) the prediction of Conjecture~\ref{conj:qh-ec} has been shown to be true
with the last step recently accomplished in~\cite{S21}.

\begin{theorem}
Let~$X$ be the adjoint variety of a non simply laced simple complex algebraic group~$\rG$.
Then $\Db(X)$ has a full $\Aut(X)$-invariant rectangular Lefschetz exceptional collection with zero residual category.
\end{theorem}

The proof is a combination of~\cite[Theorem~7.1]{K08-IGr} for type~$\rB_n$, 
\cite[Example~1.4]{KS20} for type~$\rC_n$, 
\cite[Theorem~1.1]{S21} for type~$\rF_4$, 
and~\cite[\S6.4]{K06} for type~$\rG_2$.

\subsection{Residual categories of hypersurfaces}
\label{ss:residual-hypersurfaces}

Until now we only met examples of residual categories which were of combinatorial nature 
(either generated by completely orthogonal exceptional collections or equivalent to derived categories of Dynkin quivers).
In~\S\S\ref{ss:residual-hypersurfaces}--\ref{ss:residual-ci}
we consider more complicated examples and
concentrate on a description of their basic invariants: \emph{Serre functors} and \emph{Serre dimensions}.

Let~$X \subset \P^{n}$ be a Fano hypersurface, i.e., a hypersurface of degree~$d \le n$.
Note that
\begin{equation*}
\omega_X^{-1} \cong \cO_X(n + 1 -d), 
\end{equation*}
hence 
\eqref{eq:omega-root} holds 
for~$\cL = \cO_X(1)$ and
\begin{equation*}
m = n + 1 - d.
\end{equation*}
Furthermore, the category~$\cB \coloneqq \langle \cO_X \rangle$ is admissible in~$\Dp(X)$ 
(we do not assume~$X$ to be smooth and therefore consider~$\Dp(X)$ instead of~$\Db(X)$)
and induces a rectangular Lefschetz collection of length~$m$ 
\begin{equation}
\label{eq:residual-simple-fano}
\Dp(X) = \langle \cR_X, \cO_X, \dots, \cO_X(m-1) \rangle
\end{equation}
defining the residual category~$\cR_X \subset \Dp(X)$.
Finally, as~$X$ is a Gorenstein scheme, the category~\mbox{$\cC = \Dp(X)$} has a Serre functor 
given by tensor product with~$\omega_X$ and shift by the dimension~\mbox{$\dim(X) = n - 1$}, 
so if~$\bal_\cC$ is defined by~\eqref{eq:balpha} then the autoequivalence
\begin{equation*}
\bS_\cC \circ \bal_\cC^m \cong [n-1]
\end{equation*}
preserves any triangulated subcategory~$\cB \subset \Dp(X)$;
in particular, Serre compatibility~\eqref{eq:serre-compatibility} holds for~$\cB$ as above.

One of the most surprising properties of the residual category~$\cR_X$ of a Fano hypersurface is its fractional Calabi--Yau property
(note, however, that the residual categories that appeared in~\S\ref{ss:residual-homogeneous} are also fractional Calabi--Yau).
This result has already been explained in~\cite{K14} for smooth~$X$; we restate it here for completeness.

\begin{theorem}[{\cite[Corollary~4.3]{K04}}]
\label{theorem:serre-residual-hypersurface}
Let~$X \subset \P^n$ be a smooth hypersurface of degree~\mbox{$1 \le d \le n$}.
Set~$c = \gcd(d,n+1)$. 
Then
\begin{equation}
\label{eq:serre-residual-hypersurface}
\bS_{\cR_X}^{d/c} \cong [(n+1)(d-2)/c].
\end{equation}
\end{theorem}

\begin{remark}
The result also holds true (trivially) when~$d = n + 1$.
Indeed, the definition~\eqref{eq:residual-simple-fano} in this case implies~$\cR_X = \Dp(X)$
and the formula~\eqref{eq:serre-residual-hypersurface} reads as~$\bS_{\cR_X} \cong [n - 1]$,
which is true since~$X$ is a Calabi--Yau variety of dimension~$n - 1$.
\end{remark}

Actually, the smoothness assumption in the statement of Theorem~\ref{theorem:serre-residual-hypersurface} can be removed; 
this follows immediately from Theorem~\ref{theorem:serre-residual-general} and Remark~\ref{rem:crm-zero} below.

\begin{example}
\label{ex:residual-cubic-4d}
If~$X$ is a cubic fourfold, one has~$\bS_{\cR_X} \cong [2]$. 
Thus, the category~$\cR_X$ is a \emph{K3~category}, see~\cite{K10}.
\end{example}

In~\cite[Theorem~3.5]{K19} Theorem~\ref{theorem:serre-residual-hypersurface} was generalized to the situation
where~$X$ is a smooth divisor in (or a double covering of)
a smooth variety~$M$ which admits a rectangular Lefschetz decomposition (so that the residual category of~$M$ is zero).
We do not state this result separately because it is a special case of Theorem~\ref{theorem:serre-residual-general} stated below,
see Remark~\ref{rem:crm-zero}.
The special case of Theorem~\ref{theorem:serre-residual-hypersurface} 
is obtained by taking~$M = \P^n$ with the rectangular Lefschetz decomposition given by the Beilinson exceptional collection
\begin{equation*}
\Db(\P^n) = \langle \cO, \cO(1), \dots, \cO(n) \rangle.
\end{equation*}
There are many other special cases of~\cite[Theorem~3.5]{K19} (see~\cite[\S4]{K19} for a list)
which explain most of the currently known examples of fractional Calabi--Yau categories.  
For instance, it explains the appearance of K3~categories in derived categories of cubic fourfolds, 
Gushel--Mukai varieties of even dimensions, and Debarre--Voisin 20-folds, see~\cite[\S4.4]{K19}.

\subsection{Residual categories of complete intersections}
\label{ss:residual-ci}

The results of~\cite{K19} have been significantly generalized in~\cite{KP21}.
To explain this generalization recall that for any (enhanced) functor~$\Psi \colon \cC \to \cD$ 
between (enhanced) triangulated categories with a right adjoint functor~$\Psi^!$
(for instance, for a Fourier--Mukai functor between (perfect) derived categories of projective varieties)
one can define \textsf{twist functors}~\mbox{$\bT_{\Psi,\Psi^!} \colon \cD \to \cD$} and~$\bT_{\Psi^!,\Psi} \colon \cC \to \cC$ 
by means of distinguished triangles of functors
\begin{equation*}
\Psi \circ \Psi^! \xrightarrow{\ \eps\ } \id_\cD \xrightarrow{\quad} \bT_{\Psi,\Psi^!} 
\qquad\text{and}\qquad 
\bT_{\Psi^!,\Psi} \xrightarrow{\quad} \id_\cC \xrightarrow{\ \eta\ } \Psi^! \circ \Psi,
\end{equation*}
where~$\eta$ is the unit and~$\eps$ is the counit of adjunction.
The functor~$\Psi$ is called \textsf{spherical} if the twist functors~$\bT_{\Psi,\Psi^!}$ and~$\bT_{\Psi^!,\Psi}$ are both autoequivalences
(for alternative definitions and characterizations of spherical functors see~\cite{Anno13,AL17,K19,KSS20});
in this case the twist functors are known as \textsf{spherical twists}.
Note that if the functor~$\Psi$ is zero (for instance, if the source or target category of~$\Psi$ is zero),
it is spherical and the corresponding spherical twists are isomorphic to the identity.

The simplest geometric example of a spherical functor is the pullback functor
\begin{equation*}
\Psi \coloneqq i^* \colon \Dp(M) \to \Dp(X)
\end{equation*}
for a divisorial embedding~$i \colon X \hookrightarrow M$.
Its right adjoint is the pushforward~$i_*$ and the corresponding spherical twists are given by 
\begin{equation*}
\bT_{i^*,i_*}(\cF) = \cF \otimes \cO_X(-X)[2]
\qquad\text{and}\qquad 
\bT_{i_*,i^*}(\cG) = \cG \otimes \cO_M(-X).
\end{equation*}
Another interesting example is the pullback functor for a flat double covering~$f \colon X \to M$
(see~\cite[Lemma~2.9]{KP21} for the description of the corresponding spherical twists).

Now assume~$M$ is a projective Gorenstein variety such that~$\omega_M \cong \cL_M^{-m}$ for a line bundle~$\cL_M$
and with a semiorthogonal decomposition
\begin{equation*}
\Dp(M) = \langle \cR_M, \cB_M, \cB_M \otimes \cL_M, \dots, \cB_M \otimes \cL_M^{m-1} \rangle.
\end{equation*}
Assume furthermore given another projective Gorenstein variety~$X$ and a spherical functor
\begin{equation*}
\Psi \colon \Dp(M) \to \Dp(X)
\end{equation*}
such that~$\bT_{\Psi^!,\Psi}(\cB_M) = \cB_M \otimes \cL_M^d$ for some~$1 \le d \le m - 1$ 
and~$\Psi(- \otimes \cL_M) \cong \Psi(-) \otimes \cL_X$ for a line bundle~$\cL_X$ on~$X$.
Under these assumptions the following result is proved.

\begin{theorem}[{\cite[Corollary~4.19]{KP21}}]
\label{theorem:serre-residual-general}
Assume~$M$, $\cB_M$, $\cL_M$, $X$, $\cL_X$, and~$\Psi$ are as above.
\begin{enumerate}[label=\textup{(\roman*\textup)}]
\item 
The functor~$\Psi\vert_{\cB_M}$ is fully faithful, the subcategory~$\cB_X \coloneqq \Psi(\cB_M) \subset \Dp(X)$ is admissible,
and there is a semiorthogonal decomposition
\begin{equation*}
\Dp(X) = \langle \cR_X, \cB_X, \cB_X \otimes \cL_X, \dots, \cB_X \otimes \cL_X^{m-d-1} \rangle,
\end{equation*}
where~$\cR_X \subset \Dp(X)$ is the residual category.
\item 
The restriction~$\Psi_\cR \coloneqq \Psi\vert_{\cR_M}$ is a spherical functor~$\cR_M \to \cR_X$ between the residual categories.
\item 
If~$c = \gcd(d,m)$ then
\begin{equation}
\label{eq:serre-functors}
\begin{aligned}
\bS_{\cR_M}^{d/c} &\cong \bT_{\Psi_\cR^!,\Psi_\cR}^{m/c\hphantom{(-d)}} \circ \Big[\tfrac{d\dim(M)}c\Big],\\
\bS_{\cR_X}^{d/c} &\cong \bT_{\Psi_\cR,\Psi_\cR^!}^{(m-d)/c} \circ \Big[\tfrac{d\dim(X) - 2(m-d)}c\Big],
\end{aligned}
\end{equation} 
where~$\bT_{\Psi_\cR^!,\Psi_\cR}$ and~$\bT_{\Psi_\cR,\Psi_\cR^!}$ are the spherical twists 
with respect to the spherical functor~$\Psi_\cR$.
\end{enumerate}
\end{theorem}

\begin{remark}
\label{rem:crm-zero}
In the special case, where~$\cR_M = 0$, the spherical twist~$\bT_{\Psi_\cR,\Psi_\cR^!}$ is isomorphic to the identity,
and we conclude that~$\cR_X$ is a fractional CY category of dimension~\mbox{$\dim(X) - 2(m-d)/d$}.
\end{remark}

\begin{example}
\label{ex:serre-23}
Assume the base field is algebraically closed of characteristic not equal to~$2$.
Let~$M \subset \P^5$ be a smooth quadric and let~$X \subset M$ be its smooth intersection with a cubic hypersurface.
Then the residual category of~$M$ is generated by two completely orthogonal spinor bundles of rank~$2$:
\begin{equation*}
\cR_M = \langle \cS_+, \cS_- \rangle.
\end{equation*}
Their restrictions~$\cS_{+X}$ and~$\cS_{-X}$ to~$X$ are contained in~$\cR_X$ and form a so-called \textsf{spherical pair} 
(i.e., induce a spherical functor from the derived category of a disjoint union of two points to~$\cR_X$, see~\cite[\S2.2]{KP21}),
and the formula~\eqref{eq:serre-functors} gives
\begin{equation*}
\bS_{\cR_X}^3 \cong \bT_{\cS_{+X},\cS_{-X}} \circ [7],
\end{equation*}
where~$\bT_{\cS_{+X},\cS_{-X}} \colon \cR_X \to \cR_X$ is the spherical twist with respect to the spherical pair, i.e.,
\begin{equation*}
\bT_{\cS_{+X},\cS_{-X}}(\cF) \cong 
\Cone\Big(\Ext^\bullet(\cS_{+X},\cF) \otimes \cS_{+X} \oplus \Ext^\bullet(\cS_{-X}, \cF) \otimes \cS_{-X} \xrightarrow{\quad} \cF\Big).
\end{equation*}
\end{example}

\begin{example}
\label{ex:serre-23-refined}
In the situation of Example~\ref{ex:serre-23} a similar result 
can be proved for the \textsf{refined residual category}~$\cA_X$ of~$X$ 
defined as the orthogonal of one of the spinor bundles in~$\cR_X$, say~$\cS_{+X}$, which is exceptional, 
so that there is a semiorthogonal decomposition
\begin{equation*}
\cR_X = \langle \cA_X, \cS_{+X} \rangle.
\end{equation*}
In this case the projection of the other spinor bundle to~$\cA_X$ is a spherical object~$\rK \in \cA_X$ 
and it is proved in~\cite[Proposition~5.18]{KP21} that
\begin{equation*}
\bS_{\cA_X}^3 \cong \bT_\rK^{-1} \circ [7]
\end{equation*}
(where~$\bT_\rK$ is the spherical twist with respect to~$\rK$), 
quite similarly to the case of~$\cR_X$.
\end{example}

One can apply Theorem~\ref{theorem:serre-residual-general} in order to compute Serre dimensions of residual categories of complete intersections.
Recall that the \textsf{upper and lower Serre dimensions}~$\usdim(\cC)$ and~$\lsdim(\cC)$ 
of a category~$\cC$ admitting a Serre functor~$\bS_\cC$ 
are defined as the rate of growth of upper and lower cohomological amplitude of powers of~$\bS_\cC^{-1}$,
see~\cite{EL21} or~\cite[Definition~6.10]{KP21} for details.
In the case where~$\cC = \Dp(X)$ for a Gorenstein variety~$X$, so that~\mbox{$\bS_\cC(\cF) \cong \cF \otimes \omega_X[\dim X]$}, 
one obtains from~\cite[Lemma~5.6]{EL21} the equalities 
\begin{equation*}
\usdim(\Dp(X)) = \lsdim(\Dp(X)) = \dim(X),
\end{equation*}
so Serre dimensions provide a categorical interpretation of the geometric (Krull) dimension of a variety.
In this example the upper and lower Serre dimensions coincide, but in general this is not true,
and residual categories of Fano complete intersections provide nice examples of this sort.

\begin{theorem}[{\cite[Theorem~1.7]{KP21}}]
\label{theorem:serre-dimension-ci}
Let~$X \subset \P^n$ be a smooth Fano complete intersection in~$\P^n$ of type~$(d_1,d_2,\dots,d_k)$, where 
\begin{equation*}
d_1 \ge d_2 \ge \dots \ge d_k \ge 2.
\end{equation*}
Denote by~$\ind(X) = n + 1 - \sum_{i=1}^k d_i$ the Fano index of~$X$.
Let~$\cR_X$ be the residual category of~$X$ defined by~\eqref{eq:residual-simple-fano} with~$m = \ind(X)$. 
Assume there exists a chain of smooth varieties
\begin{equation*}
X = X_k \subset \dots \subset X_2 \subset X_1 \subset X_0 = \P^n,
\end{equation*}
where~$X_i$ is a complete intersection of type~$(d_1,d_2,\dots,d_i)$.
Then
\begin{equation*}
\usdim(\cR_X) = \dim(X) - 2 \tfrac{\ind(X)}{d_1},
\qquad 
\lsdim(\cR_X) = \dim(X) - 2 \tfrac{\ind(X)}{d_k}.
\end{equation*}
In particular, if~$d_1 > d_k$, the upper Serre dimension of~$\cR_X$ is strictly bigger than the lower Serre dimension.
\end{theorem}

The assumption of the existence of a chain of smooth complete intersections~$X_i$ 
interpolating between~$\P^n$ and~$X$ is of technical nature;
it well may be that the result is also true without this assumption.
Note also that this assumption holds when the characteristic of the base field is zero by Bertini's Theorem,
see~\cite[Lemma~6.11]{KP21}.

In the situation of Examples~\ref{ex:serre-23} and~\ref{ex:serre-23-refined} 
one deduces from Theorem~\ref{theorem:serre-dimension-ci} that
\begin{equation*}
\usdim(\cA_X) = \usdim(\cR_X) = 7/3,
\qquad 
\lsdim(\cA_X) = \lsdim(\cR_X) = 2.
\end{equation*}
In fact, it is easy to descibe objects of the categories~$\cA_X$ or~$\cR_X$ 
on which the rate of growth of powers of the inverse Serre functor equals~$7/3$ and~$2$, respectively;
indeed, the first happens on the orthogonals~$\rK^\perp \subset \cA_X$ 
and~$\cS_{+X}^\perp \cap \cS_{-X}^\perp \subset \cR_X$, respectively,
while the second holds on the subcategories generated by~$\rK$ in~$\cA_X$ and~$\cS_{\pm X}$ in~$\cR_X$, respectively.


\section{Simultaneous categorical resolutions of singularities}
\label{sec:scr}

The goal of this section is to explain the proof of Theorem~\ref{thm:intro-even} 
that provides a simultaneous categorical resolution of singularities for a nodal degeneration of even-dimensional varieties.
We start by explaining what we mean by a simultaneous categorical resolution;
this notion is similar to a relative version of the definition of a categorical resolution from~\cite{K08,KL}.

Let~$f \colon \cX \to B$ be a flat proper morphism to a pointed scheme~$(B,o)$.
Recall notation~\eqref{eq:b-notation} and~\eqref{eq:diagram-family}.
We usually assume that~$B$ is a curve and~$f$ is smooth over~$B^o$.

\begin{definition}[\cite{K21}]
\label{def:scr}
A \textsf{simultaneous categorical resolution} of~$(\cX,\cX_o)$
is a triple~$(\cD,\pi^*,\pi_*)$, 
where
\begin{itemize}
 \item 
 $\cD$ is an enhanced \emph{$B$-linear} triangulated category, and
 \item 
 $\pi^* \colon \Dp(\cX) \to \cD$ and $\pi_* \colon \cD \to \Db(\cX)$ is a pair of \emph{$B$-linear} triangulated functors, 
\end{itemize}
such that 
\begin{enumerate}[label=\textup{(\roman*)}]
 \item 
 \label{item:scr-cd}
 $\cD$ is smooth and proper over~$B$,
 \item 
 \label{item:scr-adjunction}
 $\pi^*$ is left adjoint to $\pi_*$, 
 \item 
 \label{item:scr-composition}
 $\pi_* \circ \pi^* \cong \id$.
\end{enumerate}
\end{definition}

More precisely, the condition in part~\ref{item:scr-adjunction} means that there is a functorial isomorphism
\begin{equation*}
\Hom(\pi^*\cF,\cG) \cong \Hom(\cF,\pi_*\cG)
\end{equation*}
for all~$\cF \in \Dp(\cX)$, $\cG \in \cD$, and the condition in part~\ref{item:scr-composition} 
means that the composition~\mbox{$\pi_* \circ \pi^*$} is isomorphic to the canonical inclusion~$\Dp(\cX) \hookrightarrow \Db(\cX)$.
Furthermore, we usually assume that the base change~$\cD_{B^o}$ of the category~$\cD$ along the open embedding~$B^o \hookrightarrow B$
is equivalent to~$\Dp(\cX^o) = \Db(\cX^o)$ via the functors induced by~$\pi^*$ and~$\pi_*$.

If the scheme~$\cX$ has rational singularities and~$f \colon \cX \to B$ 
admits a simultaneous resolution~\mbox{$\pi \colon \tcX \to \cX$} in the geometric sense
(i.e., a resolution of singularities of~$\cX$
such that its central fiber~$\tcX_o \to \cX_o$ is a resolution of singularities of~$\cX_o$)
then~$\cD \coloneqq \Dp(\cX) = \Db(\cX)$ with the derived pullback~$\pi^*$ and pushforward~$\pi_*$ functors
is a simultaneous categorical resolution.

Geometric simultaneous resolutions of singularities exist for surface degenerations with rational double points by~\cite{Bri} (see also~\cite{Tyurina}),
but not in higher dimensions.

\subsection{General results}
\label{ss:scr-general}

In~\cite{K21} we suggest a construction of a simultaneous categorical resolution, 
analogous to the construction of a categorical resolution of a variety~$X$ given in~\cite{K08}.
Recall that the construction of~\cite{K08} assumes that~$X$ is resolved by a single blowup with exceptional divisor~$E$
and as an extra input one needs a Lefschetz decomposition of~$\Db(E)$ with respect to the conormal line bundle of~$E$.

Similarly, to construct a simultaneous categorical resolution we assume that both the total space~$\cX$ and the central fiber~$\cX_o$ of~$f$ 
are resolved by blowups with the same center~\mbox{$Z \subset \cX_o \subset \cX$}, 
such that both have smooth exceptional divisors~$E$ and~$E_o$,
and that an appropriate Lefschetz decomposition of~$\Db(E)$ is given.
The precise statement is as follows:

\begin{theorem}[{\cite[Theorem~3.11]{K21}}]
\label{theorem:simultaneous}
Let~$f \colon \cX \to B$ be a flat projective morphism to a smooth pointed curve~$(B,o)$ 
such that~$\cX \times_B B^o$ is smooth over~$B^o$ 
and let~$Z \subset \cX_o$ be a smooth closed subscheme in the central fiber.
Assume~$\cX$ has rational singularities, the blowups~\mbox{$\tcX \coloneqq \Bl_{Z}(\cX)$}, \mbox{$\tcX_o \coloneqq \Bl_{Z}(\cX_o)$},
and their exceptional divisors~$E$ and~$E_o$ are all smooth, and the central fiber of~\mbox{$\tcX \to B$} is reduced.
Let~$\pi \colon \tcX \to \cX$, $\pi_o \colon \tcX_o \to \cX_o$, $p \colon E \to Z$, $p_o \colon E_o \to Z$, 
$\eps \colon E \to \tcX$, $\eps_o \colon E_o \to \tcX_o$, and~$i_E \colon E_o \to E$ be the natural morphisms, 
shown on the diagram
\begin{equation}
\label{eq:big-diagram}
\vcenter{\xymatrix@R=3ex{
& 
E \ar@{^{(}->}[dl]_\eps \ar[dd]^(.3){p} &&
E_o \ar@{_{(}->}[ll]_{i_E} \ar@{^{(}->}[dl]_{\eps_o} \ar[dd]^(.3){p_o} 
\\
\tcX \ar[dd]^(.3){\pi} &&
\tcX_o \ar@{_{(}->}[ll] \ar[dd]^(.3){\pi_o}
\\
&
Z \ar@{^{(}->}[dl] &&
Z \ar@{=}[ll] \ar@{^{(}->}[dl] 
\\
\cX \ar[d]^f &&
\cX_o \ar@{_{(}->}[ll] \ar[d]
\\
B &&
\{o\}. \ar@{_{(}->}[ll]
}}
\end{equation}
Furthermore, assume given a~$Z$-linear left Lefschetz decomposition
\begin{align}
\label{eq:dbe-scr}
\Db(E) &= \langle 
\cA_{1-m} \otimes \cO_{E}((m-1)E), \dots, 
\cA_{-2} \otimes \cO_{E}(2E), 
\cA_{-1} \otimes \cO_{E}(E), \cA_0 \rangle
\end{align}
of~$\Db(E)$ such that~$p^*(\Db(Z)) \subset \cA_0$.
Then the category 
\begin{equation}
\label{eq:scr-cd}
\cD \coloneqq \{ \cF \in \Db(\tcX) \mid \eps^*(\cF) \in \cA_0 \} \subset \Db(\tcX)
\end{equation} 
provides a categorical resolution of~$\cX$ and there is a semiorthogonal decomposition
\begin{equation}
\label{eq:scr-db-tcx}
\Db(\tcX) = \langle 
\eps_*(\cA_{1-m} \otimes \cO_E((m-1)E)), \dots, 
\eps_*(\cA_{-2} \otimes \cO_E(2E)), 
\eps_*(\cA_{-1} \otimes \cO_E(E)), 
\cD \rangle.
\end{equation}

Moreover, if additionally we have~$\cA_{-1} = \cA_0$,
and the categories~\mbox{$\cA'_k \coloneqq i_E^*(\cA_k)$}
form a semiorthogonal decomposition
\begin{equation}
\label{eq:dbeo-scr}
\Db(E_o) = \langle 
\cA'_{1-m} \otimes \cO_{E_o}((m-1)E_o), \dots, 
\cA'_{-2} \otimes \cO_{E_o}(2E_o), 
\cA'_{-1} \otimes \cO_{E_o}(E_o) \rangle,
\end{equation} 
of~$\Db(E_o)$ then:
\begin{enumerate}[label=\textup{(\roman*\textup)}]
\item 
\label{theorem:item:cd-o}
The base change~$\cD_o$ of~$\cD$ along the embedding~$\{o\} \hookrightarrow B$ is smooth and proper over~$\kk$, 
one has
\begin{equation}
\label{eq:cd-o}
\cD_o \simeq \{ \cF \in \Db(\tcX_o) \mid \eps_o^*(\cF) \in \cA'_{-1} \},
\end{equation} 
and there is a semiorthogonal decomposition
\begin{equation*}
\Db(\tcX_o) = \langle \eps_{o*}(\cA'_{1-m} \otimes \cO_{E_o}((m-2)E_o)), \dots, 
\eps_{o*}(\cA'_{-2} \otimes \cO_{E_o}(E_o)), 
\cD_o \rangle.
\end{equation*}
\item 
\label{theorem:item:cd}
The triple $(\cD,\pi^*,\pi_*)$ is a simultaneous categorical resolution of~$(\cX,\cX_o)$; 
in particular~$\cD$ is smooth and proper over~$B$.
\end{enumerate}
\end{theorem}

The category~$\cD$ defined by~\eqref{eq:scr-cd} provides a categorical resolution of~$\cX$
and fits into the semiorthogonal decomposition~\eqref{eq:scr-db-tcx} by~\cite[Theorem~4.4 and Proposition~4.1]{K08}.
The crucial step in the proof of Theorem~\ref{theorem:simultaneous} 
is the identification~\eqref{eq:cd-o} of the base change~$\cD_o$ of the category~$\cD$,
which a priori is a subcategory of the derived category of the central fiber of the morphism~$\tcX \to B$, 
with a subcategory of~$\tcX_o$.
Note that the central fiber is a reduced, but reducible scheme --- 
the union~$\tcX_o \cup E$ of the blowup~$\tcX_o$ of~$\cX_o$ and the exceptional divisor~$E$ of~$\tcX$ with~$\tcX_o \cap E = E_o$.
The required identification is achieved in~\cite[Proposition~3.7]{K21}, 
where a more general result of this sort is established.
This proves part~\ref{theorem:item:cd-o} of the theorem.

On the other hand, the right hand side of~\eqref{eq:cd-o} is an admissible subcategory of~$\Db(\tcX_o)$ by~\cite[Proposition~4.1]{K08};
in particular, it is smooth and proper.
Therefore, the second part~\ref{theorem:item:cd} of the theorem follows from the following general and very useful result.

\begin{theorem}[{\cite[Theorem~2.10]{K21}}]
\label{theorem:sp}
Let~$g \colon \cY \to B$ be a flat proper morphism of quasiprojective schemes and let
\begin{equation*}
\Db(\cY) = \langle \cD, {}^\perp\cD \rangle
\end{equation*}
be a $B$-linear semiorthogonal decomposition with admissible components and projection functors of finite cohomological amplitude.
If for each point~$b \in B$ the category~$\cD_b$ is smooth and proper over the residue field of~$b$ 
then the category~$\cD$ is smooth and proper over~$B$.
\end{theorem}

\begin{remark}
The geometric origin of the category~$\cD$ in Theorem~\ref{theorem:sp} is important for the proof given in~\cite{K21}.
We expect a similar statement for general $B$-linear categories is not true, although we do not know any example where it fails.
\end{remark}

Theorem~\ref{theorem:sp} applies immediately to the morphism $g = f \circ \pi \colon \tcX \to B$.
Indeed, for~$b \ne o$ the fiber~$\cD_b$ is equivalent to the category~$\Db(\tcX_b) = \Db(\cX_b)$
which is smooth and proper because~\mbox{$\cX \times_B B^o$} is assumed to be smooth and proper over~$B^o$,
while the category~$\cD_o$ is smooth and proper by part~\ref{theorem:item:cd-o} of the theorem.

\subsection{Nodal singularities}
\label{ss:scr-nodal}

Theorem~\ref{theorem:simultaneous} applies easily to nodal degenerations 
of even-dimensional varieties under a mild technical assumption
that can be satisfied by a simple trick, see Remark~\ref{rem:double-cover}.
We use notation introduced in diagram~\eqref{eq:big-diagram}.

\begin{theorem}[{\cite[Theorem~3.14]{K21}}]
\label{theorem:nodal-general}
Let~$\kk$ be an algebraically closed field of characteristic not equal to~$2$.
Let~\mbox{$f \colon \cX \to B$} be a flat projective morphism of relative dimension~$2n$ to a smooth pointed curve~$(B,o)$ 
such that~$\cX \times_B B^o$ is smooth over~$B^o$.
Assume the central fiber~$\cX_o$ and the total space~$\cX$ both have an isolated ordinary double point at~$x_o$.
Let~$E$ and~$E_o$ be the exceptional divisors of the blowups~$\Bl_{x_o}(\cX)$ and~$\Bl_{x_o}(\cX_o)$.
Then~$(\cX,\cX_o)$ has a simultaneous categorical resolution of singularities~$\cD$
fitting into a semiorthogonal decomposition
\begin{equation*}
\Db(\Bl_{x_o}(\cX)) = \langle \eps_*\cO_{E}((2n-1)E), \dots, \eps_*\cO_{E}(E), \eps_*\cS_{E}, \cD \rangle,
\end{equation*}
where~$\cS_E$ is a spinor bundle on the smooth quadric~$E$;
in particular, $\cD$ is smooth and proper over~$B$
with~$\cD_b = \Db(\cX_b)$ for~$b \ne o$, and the central fiber~$\cD_o$ of~$\cD$ fits into a semiorthogonal decomposition
\begin{equation*}
\Db(\Bl_{x_o}(\cX_o)) = \langle \eps_{o*}\cO_{E_o}((2n-2)E_o), \dots, \eps_{o*}\cO_{E_o}(E_o), \cD_o \rangle.
\end{equation*}
\end{theorem}

Indeed, in this case the exceptional divisor~$E$ of the blowup~$\tcX = \Bl_{x_o}(X)$ is a smooth quadric of dimension~$2n$,
so we can take~\eqref{eq:dbe-scr} to be the left Lefschetz decomposition from Lemma~\ref{lemma:quadric-lefschetz}.
Then~\mbox{$\cA_{-1} = \cA_0$} (here it is important that~$\dim(\cX/N) = \dim(E)$ is even), 
and since~$E_o$ is a smooth quadric of dimension~$2n - 1$, 
decomposition~\eqref{eq:dbeo-scr} holds, again by Lemma~\ref{lemma:quadric-lefschetz}.
So, Theorem~\ref{theorem:simultaneous} gives the desired results.

\begin{remark}
\label{rem:double-cover}
If~$f' \colon \cX' \to B'$ is a \emph{smoothing} of a nodal variety~$X$ (i.e., $\cX'_o \cong X$ and~$\cX'$ is smooth)
applying base change with respect to a double covering~$B \to B'$ ramified over~$o$,
we obtain a morphism~$\cX \to B$ such that~$\cX_o \cong X$ and~$\cX$ has an ordinary double point at~$x_o$.
This double covering trick is quite standard, see~\cite{At58}.
\end{remark}

Note that the simultaneous categorical resolution~$\cD$ constructed in Theorem~\ref{theorem:nodal-general} depends
on a choice of one of the two spinor bundles~$\cS_E$;
however two different choices result in equivalent categorical resolutions, 
and the equivalence can be thought of as an instance of a categorical flop, see~\cite[Proposition~3.15]{K21}.

In the case of a degeneration of surfaces (i.e., for~$n = 1$), 
the category~$\cD$ constructed in Theorem~\ref{theorem:nodal-general}
is equivalent to the derived category of a small resolution of singularities of~$\cX$,
a choice of one of the two small resolutions corresponds to a choice of one of the two spinor bundles~$\cS_E$ 
on the smooth quadric surface~$E$,
and the categorical flop mentioned above reduces to the usual Atiyah flop between the small resolutions.

It would be very interesting to find generalizations of Theorem~\ref{theorem:nodal-general} to other types of simple singularities.

\subsection{Application to nodal degenerations of cubic fourfolds}
\label{ss:scr-cubic}

In this subsection we provide a simple application of Theorem~\ref{theorem:nodal-general} to K3~categories of cubic fourfolds.
Recall that for any cubic fourfold~\mbox{$X \subset \P^5$} there is a semiorthogonal decomposition 
\begin{equation}
\label{eq:sod-cubic-4}
\Dp(X) = \langle \cR_X, \cO_X, \cO_X(1), \cO_X(2) \rangle.
\end{equation}
In fact, this is just the special case of decomposition~\eqref{eq:residual-simple-fano},
so the category~$\cR_X$ above is the residual category of~$X$.
In particular, as we observed in Example~\ref{ex:residual-cubic-4d},
the Serre functor of~$\cR_X$ is isomorphic to the shift~$[2]$, so~$\cR_X$ is a \emph{K3~category}.

In the case where~$X$ has a single ordinary double point~$x_o \in X$, 
the category~$\cR_X$ is not smooth, but it admits a categorical resolution by the derived category of a smooth K3~surface.
In fact, the linear projection out of~$x_o$ identifies the blowup~$\Bl_{x_o}(X)$ 
with the blowup~$\Bl_S(\P^4)$ of~$\P^4$ along a complete intersection~K3~surface~$S \subset \P^4$ of type~$(2,3)$
and the derived category of~$S$ provides a categorical resolution of singularities for~$\cR_X$, see~\cite[Theorem~5.2]{K10}. 
We write~$S(X,x_o)$ for this K3~surface.

In the next theorem we show that the category~$\Db(S(X,x_o))$ can be realized as a limiting category
for the K3~categories of smooth cubic fourfolds.

\begin{theorem}[{\cite[Corollary~1.8]{K21}}]
\label{thm:cubic-4}
Let~$X$ be a cubic fourfold with a single ordinary double point~$x_o \in X$ 
over an algebraically closed field~$\kk$ of characteristic not equal to~$2$.
There is 
\begin{itemize}
\item 
a flat proper family~$\cX \subset \P^5_B \to B$ of cubic fourfolds over a smooth pointed curve~$(B,o)$ with central fiber~$\cX_o \cong X$
such that~$\cX$ is smooth over~$B^o$ and has an ordinary double point at~$x_o$, and
\item 
a $B$-linear category~$\cR$ smooth and proper over~$B$ such that:
\begin{enumerate}[label=\textup{(\roman*)}]
\item 
for any point~$b \ne o$ in~$B$ one has~$\cR_b \simeq \cR_{\cX_b}$, i.e., the fiber~$\cR_b$ is the~K3~category of~$\cX_b$;
\item 
one has~$\cR_o \simeq \Db(S(X,x_o))$.
\end{enumerate}
\end{itemize}
In particular, $\Db(S(X,x_o))$ is a smooth and proper extension of the family of categories~$\cR_{\cX_b}$ across the point~$o \in B$.
\end{theorem}

The construction of the family~$\cX$ is quite straightforward --- 
we take any smooth cubic fourfold~$X' \subset \P^5$ in the ambient projective space of~$X$ 
in such a way that the singular point~$x_o \in X$ does not lie on~$X'$.
Then, if~$F$ and~$F'$ are the cubic equations of~$X$ and~$X'$, 
we consider the hypersurface in~$\P^5 \times \AA^1$ given by the equation
\begin{equation*}
(1-t^2)F + t^2F' = 0,
\end{equation*}
where~$t$ is a coordinate on~$\AA^1$.
Throwing away its singular fibers (except for the fiber~$X$ over~$0 \in \AA^1$) 
we obtain the required family~$\cX \to B \subset \AA^1$.

Next, we consider the smooth and proper $B$-linear category~$\cD \subset \Db(\Bl_{x_o}(\cX))$ from Theorem~\ref{theorem:nodal-general}
and define the subcategory~$\cR \subset \cD$ by the semiorthogonal decomposition
\begin{equation}
\label{eq:sod-cd-cubic}
\cD = \langle \cR, \tf^*\Db(B), \tf^*\Db(B) \otimes \cO_{\cX/B}(1), \tf^*\Db(B) \otimes \cO_{\cX/B}(2) \rangle,
\end{equation}
where~$\tf \colon \Bl_{x_o}(\cX) \to B$ is the composition~$\Bl_{x_o}(\cX) \xrightarrow{\ \pi\ } \cX \xrightarrow{\ f\ } B$ 
of the blowup morphism and the natural projection
and~$\cO_{\cX/B}(i)$ is the pullback to~$\Bl_{x_o}(\cX)$ of the line bundle~$\cO_{\P^5_B/B}(i)$.
Applying the base change of this decomposition to various points of the base~$B$ and using~\cite[Theorem~5.2]{K10} for the point~$o$,
we obtain the required identifications of the categories~$\cR_b$ and~$\cR_o$.  

\medskip

One possible interpretation of Theorem~\ref{thm:cubic-4} is the following.
Let~$\fMcub$ be the GIT moduli space of cubic fourfolds~\cite{Laza} 
and let~$\fMcubnod \subset \fMcub$ be the divisor of singular cubic fourfolds.
Then the family of K3~categories of smooth cubic fourfolds 
(a priori defined over the open subspace~\mbox{$\fMcub \setminus \fMcubnod \subset \fMcub$}) 
extends to the general point of the boundary divisor of the root stack~$\sqrt{\fMcubnod/\fMcub}$
(the appearance of the root stack corresponds to the necessity of the double covering trick of Remark~\ref{rem:double-cover}).

\begin{remark}
It would be interesting to find analogous extensions of K3~categories of smooth cubic fourfolds 
to other special loci of the moduli space~$\fMcub$.
One particularly interesting point of~$\fMcub$ corresponds to the so-called ``chordal cubic'' (see~\cite[\S4.4]{Has} and~\cite[\S8.2]{Laza}), 
defined as the secant variety of the Veronese surface~$\upsilon_2(\P^2) \subset \P^5$.
It seems likely that to make such an extension possible it is necessary to blowup this point on the moduli space.
A general point of the exceptional divisor corresponds to a smooth sextic curve in~$\P^2$
and it is natural to expect the derived category of the double covering of~$\P^2$ 
branched at that curve to show up in the extension.
\end{remark}

\section{Absorption of singularities}
\label{sec:absorption}

The goal of this section is to introduce the notion of absorption of singularities
and to explain the proof of Theorem~\ref{thm:intro-odd}. 

\subsection{Absorption and deformation absorption}
\label{ss:def-absorption}

We start with the definition of absorption of singularities of a category.
For simplicity we restrict to the case when the category in question is~$\Db(X)$ and~$X$ is a proper (singular) variety.

\begin{definition}[\cite{KS21}]
\label{def:absorption}
We say that a subcategory~\mbox{$\cP \subset \Db(X)$} \textsf{absorbs singularities of~$X$} 
if it is admissible and the orthogonals~$\cP^\perp$ and~${}^\perp\cP$ in~$\Db(X)$ are smooth and proper.
\end{definition}

Note that, when~$\cP$ is admissible, the left and right mutation functors with respect to~$\cP$ 
induce equivalences of the orthogonals~$\cP^\perp \simeq {}^\perp\cP$,
so one of them is smooth and proper if and only if the other is so.

We think of the category~$\cP^\perp \simeq {}^\perp\cP$ as a smooth and proper ``modification'' of~$\Db(X)$.
There is a trivial example of absorption with~$\cP = \Db(X)$ and~\mbox{$\cP^\perp = {}^\perp\cP = 0$}; 
of course, this example is not interesting, and it shows that it is desirable 
to have the absorbing category~$\cP$ as small as possible (hence its orthogonals~$\cP^\perp$ and~${}^\perp\cP$ as big as possible).

The notion of absorption is ``opposite'' to the notion of categorical resolution 
in the sense that in the latter we replace~$\Db(X)$ by a \emph{larger} smooth and proper category,
while in the former we replace it by a \emph{smaller} smooth and proper category.

The following is the simplest example of absorption.

\begin{example}
\label{ex:minus-one-curve}
Let $X = X_1 \cup X_2$ be a complete curve with two smooth components 
intersecting transversely at a point~$x_0$ and with~\mbox{$X_1 \cong \P^1$}.
Then
\begin{equation*}
\cP \coloneqq \langle \cO_{X_1}(-1) \rangle \subset \Db(X)
\end{equation*}
absorbs singularities of~$X$; indeed, the category~${}^\perp\cP$ is equivalent to~$\Db(X_2)$
(via the pullback functor for the projection~$X \to X_2$ contracting~$X_1$ to the intersection point~$x_0 \in X_2$),
hence smooth, and~$\cP^\perp \simeq {}^\perp\cP$ because~$X$ is Gorenstein.
\end{example}

Other examples of absorption are given by the so-called \emph{Kawamata type} semiorthogonal decompositions 
introduced in~\cite[Definition~4.1]{KPS20} (see~\cite{KKS20} for many decompositions of this type for surfaces).

We will give more examples of absorption in the next subsection, and meanwhile we introduce a stronger notion.
Recall that a \textsf{smoothing} of a proper variety~$X$ is a Cartesian diagram~\eqref{eq:diagram-family} 
where~$f$ is a flat proper morphism to a smooth pointed curve~$(B,o)$ and~$\cX$ is smooth.
Recall that~$i \colon X \to \cX$ denotes the embedding of central fiber.

\begin{definition}[\cite{KS21}]
\label{def:deformation-absorption}
Assume a subcategory~$\cP \subset \Db(X)$ absorbs singularities of a proper variety~$X$.
We say that~$\cP$ \textsf{provides a deformation absorption of singularities of~$X$} if for any smoothing~\mbox{$f \colon \cX \to B$} of~$X$
the triangulated subcategory~$\langle i_*(\cP) \rangle \subset \Db(\cX)$ 
generated in~$\Db(\cX)$ by the image of~$i_* \colon \cP \to \Db(\cX)$ is admissible.
\end{definition}

\begin{example}
\label{ex:minus-one-curve-smoothing}
The absorption of singularities described in Example~\ref{ex:minus-one-curve} is a deformation absorption, 
because for any smoothing~$\cX \to B$ of the reducible curve~$X = X_1 \cup X_2$ we have
\begin{equation*}
\cN_{X_1/\cX} \cong \cO_{X_1}(X_1) \cong \cO_{X_1}(-X_2) \cong \cO_{X_1}(-x_0) \cong \cO_{X_1}(-1),
\end{equation*}
hence~$X_1 \subset \cX$ is a $(-1)$-curve on a smooth surface~$\cX$, 
so~$i_*(\cO_{X_1}(-1)) \in \Db(\cX)$ is an exceptional object,
and hence the subcategory~$\langle i_*(\cP) \rangle = \langle i_*\cO_{X_1}(-1) \rangle$ is admissible.
\end{example}

The following result demonstrates how a subcategory providing a deformation absorption can be used
to construct a smooth family of categories.

\begin{theorem}[{\cite{KS21}}]
\label{thm:deformation-absorption}
Assume a subcategory~$\cP \subset \Db(X)$ provides a deformation absorption of singularities of a proper variety~$X$.
Let~$f \colon \cX \to B$ be a smoothing of~$X$ with~$\cX$ quasiprojective.
Define the subcategory~$\cD \subset \Db(\cX)$ from the semiorthogonal decomposition
\begin{equation}
\label{eq:cd-absorption}
\Db(\cX) = \langle \langle i_*(\cP) \rangle, \cD \rangle.
\end{equation}
Then~$\cD_b = \Db(\cX_b)$ for~$b \ne o$, and~$\cD_o = {}^\perp\cP \subset \Db(X)$.
In particular, $\cD$ is smooth and proper over~$B$.
\end{theorem}

The main thing to prove here is the equality
\begin{equation*}
\langle i_*(\cP) \rangle_o = \cP
\end{equation*}
of the base change of the category~$\langle i_*(\cP) \rangle$ with~$\cP$ (as subcategories of~$\Db(X)$).
One inclusion follows from~\cite[Corollary~5.7]{K11}.
To prove the other inclusion consider the distinguished triangle
\begin{equation}
\label{eq:rrr-triangle}
i^*i_*(\cF) \to \cF \to \cF[2]
\end{equation} 
that exists for any~$\cF \in \cP$. 
It implies that~$i^*(\langle i_*(\cP) \rangle) \subset \cP$
and using the definition of base change and admissibility of~$\cP$ the required inclusion easily follows.

\begin{example}
In the setup of Example~\ref{ex:minus-one-curve-smoothing} if~$\bcX$ is the surface 
obtained from a smoothing~$\cX$ of~$X$ by contracting the $(-1)$-curve~$X_1$ 
then~$\bcX$ is smooth, $\cX$ is isomorphic to the blowup of~$\bcX$ at a point, 
and the orthogonal~$\cD$ of~$\langle i_*(\cP) \rangle$ in~$\Db(\cX)$ is equivalent to~$\Db(\bcX)$.
\end{example}

Following this example it is suggestive to think of the category~$\cD$ from Theorem~\ref{thm:deformation-absorption} 
as a categorical contraction of~$\Db(\cX)$; this point of view is developed in~\cite{KS21}.

\subsection{$\CP^\infty$-objects}
\label{ss:cp-infty}

In this section we introduce a class of objects that can be used to construct deformation absorptions of singularities.

\begin{definition}[{\cite{KS21}}]
We say that $\rP \in \Db(X)$ is a \textsf{$\CP^\infty$-object} if
\begin{equation*}
\Ext^\bullet(\rP,\rP) \cong \kk[t],
\qquad \text{where} \qquad
\deg(t) = 2.
\end{equation*}
In other words, the derived endomorphism algebra of~$\rP$ is isomorphic 
to the cohomology algebra of the topological space~$\CP^\infty$.
\end{definition}

\begin{remark}
\label{red:pinfty-q}
In~\cite{KS21} we define a more general notion of~$\P^{\infty,q}$-objects for an arbitrary positive integer~$q$
by assuming that~$\Ext^\bullet(\rP,\rP) \cong \kk[t]$ with~$\deg(t) = q$.
Such objects can be also used to absorb singularities 
(and there are many geometrically meaningful examples of such absorptions for~$q = 1$, see Remark~\ref{rem:higher-and-even}),
but they never provide deformation absorptions unless~$q = 2$.
\end{remark}

For each~$\CP^\infty$-object we define the \textsf{canonical self-extension}~$\rM$ of~$\rP$ 
from the canonical distinguished triangle 
\begin{equation}
\label{eq:mmm}
\rM  \xrightarrow{\quad} \rP \xrightarrow{\ t\ } \rP[2],
\end{equation}
where the second arrow is given by the generator~$t \in \Ext^2(\rP, \rP)$ of~$\Ext^\bullet(\rP, \rP)$.

The following characterization is useful and easy to prove.

\begin{lemma}[{\cite{KS21}}]
\label{lemma:mr-pinfty}
If~$\rP \in \Db(X)$ is a $\CP^\infty$-object and~$\rM$ is its canonical self-extension then
\begin{equation}
\label{eq:ext-m-r}
\Ext^\bullet(\rM,\rP) \cong \kk.
\end{equation}
Conversely, if objects~$\rP,\rM \in \Db(X)$ satisfy~\eqref{eq:ext-m-r} and~$\Cone(\rM \to \rP) \cong \rP[2]$,
then~$\rP$ is a~$\CP^\infty$-object.
\end{lemma}

\begin{example}
\label{ex:minus-one-curve-pinfty}
Consider the situation described in Example~\ref{ex:minus-one-curve}.
Let~$\cL_0$ be the line bundle on~$X$ that restricts to~$X_1$ as~$\cO_{X_1}(-1)$ and to~$X_2$ as~$\cO_{X_2}$.
Similarly, let~$\cL_1$ be the line bundle on~$X$ that restricts to~$X_2$ as~$\cO_{X_2}(-x_0)$ and to~$X_1$ as~$\cO_{X_1}$.
Then there is an exact sequence
\begin{equation*}
0 \to \cO_{X_1}(-1) \to \cL_1 \to \cL_0 \to \cO_{X_1}(-1) \to 0,
\end{equation*}
where the middle arrow is defined as the composition~$\cL_1 \twoheadrightarrow \cO_{X_2}(-x_0) \hookrightarrow \cL_0$.
It follows that the objects~$\rP \coloneqq \cO_{X_1}(-1)$ and~$\rM \coloneqq \Cone(\cL_1 \to \cL_0)$ 
fit into a distinguished triangle of the form~\eqref{eq:mmm}.
It is also easy to check that~\eqref{eq:ext-m-r} holds; 
therefore~$\rP$ is a $\CP^\infty$-object.
\end{example}

As the next proposition shows, $\CP^\infty$-objects induce semiorthogonal decompositions;
it is instructive to compare this with the well-known notion of~$\CP^n$-objects (see~\cite{HT06}),
which rather give autoequivalences of derived categories.

\begin{proposition}[{\cite{KS21}}]
\label{prop:pinfty-admissible}
If~$X$ is a proper Gorenstein scheme, $\rP \in \Db(X)$ is a $\CP^\infty$-object, 
and the canonical self-extension~$\rM$ of~$\rP$ is perfect,
then the subcategories
\begin{equation*}
\cP \coloneqq \langle \rP \rangle \subset \Db(X)
\qquad\text{and}\qquad 
\cM \coloneqq \langle \rM \rangle \subset \Dp(X)
\end{equation*}
are admissible in~$\Db(X)$ and~$\Dp(X)$, respectively.
Moreover,
\begin{equation*}
\cP \simeq \Db(\kk[\epsilon]/\epsilon^2)
\qquad\text{and}\qquad 
\cM \simeq \Dp(\kk[\epsilon]/\epsilon^2), 
\qquad \text{where}\quad
\deg(\epsilon) = -1,
\end{equation*}
and~$\cP \cap \Dp(X) = \cM$.
\end{proposition}

For instance, since in the situation of Example~\ref{ex:minus-one-curve-pinfty} 
the curve~$X$ is Gorenstein and the object~$\rM$ is perfect, 
we conclude that the category~$\cP \subset \Db(X)$ generated by~$\cO_{X_1}(-1)$ is admissible.
A similar computation shows that if~$X$ is a tree of rational curves, 
there is a semiorthogonal collection of~$\CP^\infty$-objects in~$\Db(X)$ absorbing its singularities.

A useful property of~$\CP^\infty$-objects is given by the following

\begin{theorem}[{\cite{KS21}}]
\label{thm:pinfty-defabs}
Let~$\rP_1,\dots,\rP_n$ be a semiorthogonal collection of $\CP^\infty$-objects in~$\Db(X)$. 
If the category~$\cP = \langle \rP_1, \dots, \rP_n \rangle$ absorbs singularities of~$X$,
it provides a deformation absorption of singularities.
\end{theorem}

To prove the theorem we consider a smoothing~$f \colon \cX \to B$ of~$X$ and check that for each~$1 \le j \le n$ the object
\begin{equation*}
\rM_j \coloneqq i^*i_*(\rP_j)
\end{equation*}
is the canonical self-extension of~$\rP_j$ 
(this follows from the triangle~\eqref{eq:rrr-triangle} for~$\cF = \rP_j$).
Then using the adjunction isomorphisms
\begin{equation*}
\Ext^\bullet(i_*\rP_j, i_*\rP_k) \cong \Ext^\bullet(\rM_j, \rP_k),
\end{equation*}
Lemma~\ref{lemma:mr-pinfty}, and the triangles~\eqref{eq:mmm}, we deduce exceptionality and semiorthogonality of~$i_*\rP_j$,
which implies admissibility of the subcategory of~$\Db(\cX)$ generated by~$i_*\rP_1,\dots,i_*\rP_n$.

Combining Theorem~\ref{thm:pinfty-defabs} with Theorem~\ref{thm:deformation-absorption} we obtain

\begin{corollary}
Let~$\rP_1,\dots,\rP_n$ be a semiorthogonal collection of $\CP^\infty$-objects in~$\Db(X)$. 
If the category~$\cP = \langle \rP_1, \dots, \rP_n \rangle$ absorbs singularities of~$X$,
then for any smoothing~$f \colon \cX \to B$ of~$X$ 
there is a semiorthogonal decomposition~
\begin{equation}
\label{eq:cpinfty-collection}
\Db(\cX) = \langle i_*\rP_1, \dots, i_*\rP_n, \cD \rangle 
\end{equation}
and the subcategory~$\cD$ defined by~\eqref{eq:cpinfty-collection}
is smooth and proper over~$B$
with~$\cD_b = \Db(\cX_b)$ for~$b \ne o$ and~$\cD_o = {}^\perp\cP \subset \Db(X)$.
\end{corollary}

\subsection{Absorption of nodal singularities}
\label{ss:absorption-nodal}

In this section we show how to construct collections of $\CP^\infty$-objects absorbing singularities of nodal varieties of odd dimension.
We concentrate on the case of threefolds, because this case is technically simpler and the main assumption has clearer geometric meaning.

\begin{definition}[{\cite{KPS20}}]
\label{def:mnf}
A threefold~$X$ with isolated singularities is called \textsf{maximally non-factorial} if the natural morphism
from the class group of Weil divisors on~$X$ to the sum of local class groups over all singular points of~$X$
\begin{equation*}
\Cl(X) \xrightarrow{\quad} \bigoplus_{x \in \Sing(X)} \Cl(X,x)
\end{equation*}
is surjective.
\end{definition}

For simplicity consider the case where~$X$ has a single ordinary double point~$x_0 \in X$.
In this case~$\Cl(X,x_0) \cong \ZZ$,
and~$X$ is non-factorial if and only if the morphism of the class groups~$\Cl(X) \to \Cl(X,x_0)$ is surjective onto a subgroup of finite index.
Thus, maximal non-factoriality is a strengthening of the usual non-factoriality property.

Further, if~$x_0 \in X$ is a non-factorial ordinary double point of a threefold, 
there exists a small resolution~$\pi \colon \tX \to X$,
i.e., a smooth threefold~$\tX$ with a projective morphism~$\pi$ such that its exceptional locus
\begin{equation}
\label{eq:l0}
L_0 \coloneqq \pi^{-1}(x_0) \cong \P^1
\end{equation}
is a smooth rational curve. 
In these terms maximal non-factoriality is equivalent to the existence of a line bundle~$\cL$ on~$X$ such that
\begin{equation}
\label{eq:cl-l0}
\cL\vert_{L_0} \cong \cO_{L_0}(-1).
\end{equation}

\begin{theorem}[{\cite{KS21}}]
\label{thm:nmnf3}
Let~$X$ be a maximally non-factorial proper threefold with a single ordinary double point~$x_0 \in X$ and assume~$\rH^\bullet(X,\cO_X) = \kk$.
Let~$\pi \colon \tX \to X$ be a small resolution with exceptional locus~\eqref{eq:l0}.
Then for any line bundle~$\cL$ on~$\tX$ for which~\eqref{eq:cl-l0} holds the object
\begin{equation}
\label{eq:pinfty-piscl}
\rP \coloneqq \pi_*\cL
\end{equation}
is a $\CP^\infty$-object providing a deformation absorption of singularities of~$X$.
\end{theorem}

To prove the theorem we first decompose the derived category~$\Db(\tX)$ into two parts:
the first is generated by an exceptional pair and the second is its orthogonal.
The pair consists of the line bundle~$\cL$ and the twisted ideal sheaf
\begin{equation*}
\cL' \coloneqq 
\cI_{L_0} \otimes \cL.
\end{equation*}
Exceptionality of the pair~$(\cL',\cL)$ follows from the isomorphism~\eqref{eq:cl-l0} and the fact that~$L_0$ is a $(-1,-1)$-curve.
Now we obtain a semiorthogonal decomposition
\begin{equation*}
\Db(\tX) = \langle \tcP, \tcD \rangle,
\qquad\text{where}\qquad 
\tcP \coloneqq \langle \cL',\cL \rangle
\qquad\text{and}\qquad 
\tcD \coloneqq {}^\perp\langle \cL',\cL \rangle.
\end{equation*}
Since by~\cite[Theorem~2.14]{BKS18} the functor~$\pi_* \colon \Db(\tX) \to \Db(X)$ is a Verdier localization 
with the kernel generated by~$\cO_{L_0}(-1)$ and since~$\cO_{L_0}(-1) \cong \Cone(\cL' \to \cL) \in \tcP$, 
it follows that~\mbox{$\tcD = \pi^*(\cD)$} and there is a semiorthogonal decomposition
\begin{equation*}
\Db(X) = \langle \cP, \cD \rangle,
\end{equation*}
where $\cP \simeq \tcP / \langle \cO_{L_0}(-1) \rangle$.
The category~$\cD$ is equivalent to the smooth and proper category~$\tcD$, hence~$\cP$ absorbs singularities of~$X$.
On the other hand, $\cP$ is generated by the object~\mbox{$\rP \coloneqq \pi_*(\cL) \cong \pi_*(\cL')$}, 
and a simple computation shows that there is a distinguished triangle
\begin{equation*}
\pi^*(\rM) \to \cL' \to \cL[2]
\end{equation*}
for an object~$\rM \in \Dp(X)$ such that~\eqref{eq:ext-m-r} holds.
Pushing forward this triangle we obtain~\eqref{eq:mmm}, hence~$\rP$ is a $\CP^\infty$-object by Lemma~\ref{lemma:mr-pinfty}
(and~$\rM$ is its canonical self-extension).

\begin{remark}
The~$\CP^\infty$-object~$\rP$ defined in~\eqref{eq:pinfty-piscl} is a reflexive sheaf of rank~$1$ on~$X$
and the image of~\mbox{$[\rP] \in \Cl(X)$} in~$\Cl(X,x_0)$ is a generator of the local class group.
Note that the line bundle~$\cL$ satisfying~\eqref{eq:cl-l0} is unique up to twist by a line bundle pulled back from~$X$,
hence the same is true for the reflexive sheaf~$\rP$.
On the other hand, if~$\tX'$ is the other small resolution of~$X$ and~$\rP'$ 
is the reflexive generator of~$\Cl(X, x_0)$ constructed from it,
it follows that~$\rP'$ is isomorphic to the dual~$\rP^\vee$ of~$\rP$ up to line bundle twist.
Thus, the flop~$\tX \dashrightarrow \tX'$ between the small resolutions
corresponds to dualization of the corresponding~$\CP^\infty$-object.
\end{remark}

\begin{remark}
Theorem~\ref{thm:nmnf3} shows that maximal non-factoriality is sufficient 
for the existence of an absorption of singularities of a threefold with a single ordinary double point 
by a~$\CP^\infty$-object.
Using~\cite[Lemma~1.11]{PS18} and~\cite[Corollary~3.8]{KPS20} one can prove that this condition is also necessary, see~\cite{KS21}.
\end{remark}

\begin{remark}
\label{rem:mnf-many}
Let~$X$ be a maximally non-factorial threefold with~$n$ ordinary double points 
and let~$\pi \colon \tX \to X$ be a small resolution of singularities with exceptional curves~$L_1,\dots,L_n$.
Assuming there is an exceptional collection~$(\cL_1,\dots,\cL_n)$ of line bundles on the blowup~$\tX$ 
such that~\mbox{$\cL_j\vert_{L_k} \cong \cO_{L_k}(-\delta_{jk})$} for all~$1 \le j,k \le n$,
a similar argument shows that singularities of~$X$ are absorbed 
by the semiorthogonal collection of~$\CP^\infty$-objects~$\rP_j \coloneqq \pi_*\cL_j$, $1 \le j \le n$.
\end{remark}

Combining Theorem~\ref{thm:nmnf3} with Theorem~\ref{thm:pinfty-defabs} we obtain

\begin{corollary}
\label{cor:nmnf3}
Let~$X$ be a maximally non-factorial threefold with a single ordinary double point~$x_0 \in X$ and assume~$\rH^\bullet(X,\cO_X) = \kk$.
Let~$f \colon \cX \to B$ be a smoothing of~$X$.
Then there is a semiorthogonal decomposition
\begin{equation*}
\Db(\cX) = \langle i_*\rP, \cD \rangle,
\end{equation*}
where~$\rP \in \Db(X)$ is the $\CP^\infty$-object defined in Theorem~\textup{\ref{thm:nmnf3}}
and the subcategory~$\cD$ is smooth and proper over~$B$ 
with~$\cD_b = \Db(\cX_b)$ for~$b \ne o$ and~$\cD_o = {}^\perp\rP \subset \Db(X)$.
\end{corollary}

Of course, a similar result holds for maximally non-factorial threefolds with several ordinary double points 
if the condition of Remark~\ref{rem:mnf-many} is satisfied.

\begin{remark}
\label{rem:higher-and-even}
There is a similar construction of deformation absorption that works in higher dimensions.
Let~$X$ be a variety of \emph{odd} dimension with a single ordinary double point~$x_0 \in X$.
Let~$\tX = \Bl_{x_0}(X)$; then the exceptional divisor~$E \subset \tX$ is a smooth even-dimensional quadric.
Assume there is an exceptional object~$\cE \in \Db(\tX)$ such that
\begin{equation*}
\cE\vert_E \cong \cS_E,
\end{equation*}
where~$\cS_E$ is a spinor bundle (this condition plays the same role as~\eqref{eq:cl-l0}).
Then~$\rP \coloneqq \pi_*(\cE)$ is a $\CP^\infty$-object providing a deformation absorption of singularities of~$X$, see~\cite{KS21}.
Of course, there is also a version of this result for several ordinary double points as in Remark~\ref{rem:mnf-many}.

An analogous construction for \emph{even-dimensional} varieties produces a~$\P^{\infty,1}$-object (as defined in Remark~\ref{red:pinfty-q})
which also absorbs singularities of~$X$, but it does not give a deformation absorption.
\end{remark}

\subsection{Fano threefolds}
\label{ss:fano}

In this subsection we apply the above results to clarify and extend the relation 
between nontrivial components of derived categories of del Pezzo threefolds and prime Fano threefolds
that was discovered in~\cite{K09}.

Recall that a \textsf{prime Fano threefold} is a Fano threefold~$X$ with~$\Pic(X) = \ZZ K_X$,
and its \textsf{genus}~$\g(X)$ is defined from the equality
\begin{equation*}
(-K_X)^3 = 2\g(X) - 2.
\end{equation*}
It is well known that~$2 \le \g(X) \le 12$ and~$\g(X) \ne 11$.

Mukai proved (\cite{Mukai}, \cite[\S B.1]{KPS18})
that for any prime Fano threefold~$X$ with
\begin{equation*}
\g(X) \in \{6,8,10,12\}
\end{equation*}
there exists a unique exceptional vector bundle~$\cU_X$ of rank~2 with~$\rc_1(\cU_X) = K_X$
such that~$(\cO_X,\cU_X^\vee)$ is an exceptional pair; it is called the \textsf{Mukai bundle}.
Using this, the \textsf{nontrivial part}~\mbox{$\cA_X \subset \Db(X)$} 
was defined in~\cite{K09} from the semiorthogonal decomposition
\begin{equation}
\label{eq:sod-fano-3-even-genus}
\Db(X) = \langle \cA_X, \cO_X, \cU_X^\vee \rangle.
\end{equation} 

\begin{remark}
\label{rem:genus-4}
This definition extends to \emph{general} prime Fano threefolds of genus
\begin{equation*}
\g(X) = 4.
\end{equation*}
In fact, any smooth prime Fano of genus~$4$ is a complete intersection~$X \subset \P^5$ of type~$(2,3)$.
We will say that~$X$ is \textsf{general} if the (unique) quadric passing through~$X \subset \P^5$ is smooth;
in this case there are two Mukai bundles (the restrictions of the spinor bundles from the quadric) 
and the corresponding nontrivial parts of~$\Db(X)$ are equal to the refined residual categories from Example~\ref{ex:serre-23-refined}.
\end{remark}

Similarly, a \textsf{del Pezzo threefold} is a Fano threefold~$Y$ with~$-K_Y = 2H$ for a primitive Cartier divisor class~$H$
and its \textsf{degree}~$\dg(Y)$ is defined as
\begin{equation*}
\dg(Y) = H^3.
\end{equation*}
It is well known that~$1 \le \dg(Y) \le 5$ for del Pezzo threefolds of Picard rank~$1$.

If~$Y$ is a del Pezzo threefold, the pair of line bundles~$(\cO_Y,\cO_Y(H))$ is exceptional, 
and this time the \textsf{nontrivial part}~$\cB_Y \subset \Db(Y)$ was defined in~\cite{K09} from the semiorthogonal decomposition
\begin{equation}
\label{eq:sod-dp-3}
\Db(Y) = \langle \cB_Y, \cO_Y, \cO_Y(H) \rangle
\end{equation} 
(so, in this case these are just the residual categories in the sense of~\S\ref{sec:residual}).

It was observed in~\cite[Proposition~3.9]{K09} that when~$X$ and~$Y$ are as above and
\begin{equation*}
\g(X) = 2\dg(Y) + 2,
\end{equation*}
the categories~$\cA_X$ and~$\cB_Y$ have isomorphic numerical Grothendieck groups 
(and their isomorphism is compatible with the Euler pairings).
Furthermore, it was proved in~\cite[Theorem~3.8]{K09} that for each prime Fano threefold~$X$ with~$\g(X) \in \{8,10,12\}$
there is a unique del Pezzo threefold~$Y$ (with~\mbox{$\dg(Y) = \g(X)/2 - 1 \in \{3,4,5\}$}) such that
\begin{equation*}
\cA_X \simeq \cB_Y.
\end{equation*}
So, it was expected~\cite[Conjecture~3.7]{K09} that the same equivalence takes place 
for appropriate pairs~$(X,Y)$ with~$\g(X) \in \{4,6\}$ and~$\dg(Y) \in \{1,2\}$.

However, the conjecture turned out to be false:
for~$\g(X) = 6$ and~$\dg(Y) = 2$ it was disproved in~\cite{BP} or~\cite[Theorem~1.2]{Zhang}, 
and for~$\g(X) = 4$ and~$\dg(Y) = 1$ it is false for trivial reasons 
as in this case~$\dim(\HOH_1(\cA_X)) = 21$ while~$\dim(\HOH_1(\cB_Y)) = 20$.
The next theorem clarifies the situation in these two cases.

\begin{theorem}[\cite{KS21:Fano}]
\label{thm:fano-3}
For~$2 \le d \le 5$ let~$Y$ be a del Pezzo threefold of degree~$d$ and Picard rank~$1$
and for~\mbox{$d = 1$} let~$Y$ be a small resolution of a del Pezzo threefold of degree~$1$ and Picard rank~$1$
with a single ordinary double point.
Then there exists a flat projective morphism~$f \colon \cX \to B$ to a smooth pointed curve~$(B,o)$ 
such that~$\cX$ is smooth and
\begin{enumerate}[label=\textup(\alph*\textup)]
\item 
for any point~$b \ne o$ in~$B$ the fiber~$\cX_b$ is a smooth prime Fano threefold of genus~$g = 2d + 2$;
\item 
the central fiber~$\cX_o$ is a prime Fano threefold of genus~$g = 2d + 2$ birational to~$Y$
with a single maximally non-factorial ordinary double point~$x_o \in \cX_o$. 
\end{enumerate}
Furthermore, there is a $B$-linear subcategory~$\cA \subset \Db(\cX)$ which is smooth and proper over~$B$ and such that:
\begin{enumerate}[label=\textup{(\roman*\textup)}]
\item 
for any point~$b \ne o$ in~$B$ one has~$\cA_b = \cA_{\cX_b} \subset \Db(\cX_b)$;
\item 
the central fiber~$\cA_o$ is equivalent to the component~$\cB_Y$ of~$\Db(Y)$.
\end{enumerate}
In particular, the nontrivial part~$\cB_Y$ of~$\Db(Y)$ is a smooth and proper extension across the point~$o \in B$
of the family~$\cA_{\cX_b}$ of the nontrivial parts of~$\Db(\cX_b)$.
\end{theorem}

\begin{remark}
The case~$d = 1$ in the theorem is somewhat special;
in this case we take~$Y$ to be a small resolution~$Y \to \bar{Y}$ of a del Pezzo threefold~$\bar{Y}$ with a single node
(such resolution~$Y$ always exists as an algebraic space, but not as a projective variety).
Note, however, that the component~$\cB_Y$ is still defined for this algebraic space~$Y$ by the same formula~\eqref{eq:sod-dp-3},
where the line bundle~$\cO_Y(H)$ is the effective generator of the Picard group of~$Y$.
It is remarkable that in this case all prime Fano threefolds~$\cX_b$ for~$b \ne o$ in the constructed family 
are general in the sense of Remark~\ref{rem:genus-4}.
\end{remark}

Let us explain how the family~$\cX$ is constructed.
Let~$Y$ be as in the theorem.
If~\mbox{$2 \le d \le 5$} let~$C \subset Y$ be a general smooth rational curve of degree~$d - 1$ (with respect to~$H$),
and if~$d = 1$ let~$C$ be the exceptional curve of~$Y$ 
(recall that in this case~$Y$ is a small resolution of a nodal del Pezzo threefold~$\bar{Y}$).
Then one can prove that~$C$ has a unique bisecant line~$L \subset Y$ 
and there exists a diagram
\begin{equation*}
\xymatrix{
&& 
\Bl_C(Y) \ar@{=}[r] \ar[dl]_\rho &
\tX \ar[dr]^\pi &
L_0 \ar@{_{(}->}[l] \ar[dr]
\\
L \ar@{^{(}->}[r] &
Y &&&
X &
\{x_0\} \ar@{_{(}->}[l] 
}
\end{equation*}
where~$\rho$ is the blowup morphism,
$\pi$ is the contraction of the strict transform~$L_0 \subset \tX$ of~$L$,
and~$X$ is a maximally non-factorial prime Fano threefold of genus~$g = 2d + 2$
with a single ordinary double point~$x_0 = \pi(L_0)$.
Now we define the family~$f \colon \cX \to B$ as a smoothing of~$X$ (it exists by~\cite{Na97}). 

Now we explain how the subcategory~$\cA \subset \Db(\cX)$ is constructed.
First, applying Corollary~\ref{cor:nmnf3} to the smoothing~$f \colon \cX \to B$ constructed above 
we obtain a $B$-linear subcategory~$\cD \subset \Db(\cX)$ smooth and proper over~$B$ 
such that
\begin{equation*}
\cD_b = \Db(\cX_b)
\quad\text{for~$b \ne o$}
\qquad\text{and}\qquad
\cD_o \simeq {}^\perp\langle \cI_{L_0}(-H), \cO_\tX(-H) \rangle \subset \Db(\tX),
\end{equation*}
where~$H$ is the pullback to~$\tX$ of the hyperplane class of~$Y$.
Next, we consider the sheaf
\begin{equation*}
\cU'_Y \coloneqq \Ker(\cO_Y(H) \oplus \cO_Y(H) \to \cO_C(d)),
\end{equation*}
where the morphism is a twist of the evaluation morphism~$\cO_Y \oplus \cO_Y \to \cO_C(1)$.
We check that there is a vector bundle~$\cU_X$ on~$X$ such that~$\cU'_Y \cong \rho_*(\pi^*\cU_X^\vee)$,
that~$(\cO_X,\cU_X^\vee)$ is an exceptional pair in~$\Db(X)$, and 
this pair deforms (possibly after an \'etale base change) to the nearby fibers of a family~$f \colon \cX \to B$.
Therefore, after a possible \'etale base change we can assume that the pair is defined on~$\cX$,
and hence we have a $B$-linear semiorthogonal decomposition
\begin{equation*}
\cD = \langle \cA, f^*\Db(B), f^*\Db(B) \otimes \cU_\cX^\vee \rangle
\end{equation*}
Finally, we check that~$\cU_\cX\vert_{\cX_b}$ is the Mukai bundle of~$\cX_b$ when~$b \ne o$, hence~$\cA_b = \cA_{X_b}$.
On the other hand, we find a simple sequence of mutations identifying~$\cA_o \subset \cD_o \subset \Db(\tX)$ with~$\cB_Y$.
This proves an equivalence~$\cA_o \cong \cB_Y$.

\begin{remark}
When~$3 \le d \le 5$ the family of threefolds~$\cX$ in Theorem~\ref{thm:fano-3} can be chosen in such a way 
that the family of categories~$\cA_{\cX_b}$ is \emph{isotrivial},
i.e., $\cA_{\cX_b} \simeq \cB_Y$ for all~$b \in B$.
This is no longer possible for~$d \in \{1,2\}$.
\end{remark}

\begin{remark}
There are several interesting examples of del Pezzo threefolds with higher Picard rank:
two del Pezzo threefolds of degree~$6$ (the flag variety~$\Fl(1,2;3)$ and~$(\P^1)^3$)
and one del Pezzo threefold of degree~$7$ (the blowup of~$\P^3$ at a point).
The construction of Theorem~\ref{thm:fano-3} works for these threefolds 
and relates the nontrivial parts of their derived categories (still defined by~\eqref{eq:sod-dp-3})
to the nontrivial parts (still defined by~\eqref{eq:sod-fano-3-even-genus}) of the derived categories of appropriate Fano threefolds 
with primitive canonical class and genus~$g = 14$ and~$g = 16$, respectively.
\end{remark}

There are other interesting maximally non-factorial nodal Fano threefolds,
for instance some prime Fano threefolds of odd genus~\mbox{$g \in \{5,7,9\}$}.
They also provide geometrically meaningful extensions of (appropriately defined) nontrivial components 
of derived categories, see~\cite{KS21:Fano} for details.

\bibliography{./icm-bib.bib}
\bibliographystyle{plain}

\end{document}